\documentclass[11pt]{article}

\usepackage{hyperref}

\usepackage{amsmath,amsthm}
\usepackage{amssymb,bbm}
\usepackage{mathabx}
 \usepackage{IEEEtrantools}
\usepackage{mathrsfs}
\usepackage{color}
\usepackage{tikz-cd}
\usetikzlibrary{arrows,decorations.markings}
 \usetikzlibrary{fit,matrix}
\usepackage{multirow}

\usepackage{authblk}

\numberwithin{equation}{section}

\counterwithin*{equation}{section}

\newcommand{\one}{\mathbbm{1}}

\newcommand{\uqisltc}{\mathcal{U}_{q_{i}}(s\ell(2,\mathbb{C}))}

\newcommand{\utisltc}{\mathcal{U}_{t_{i}}(s\ell(2,\mathbb{C}))}

\newcommand{\utisltca}{\mathcal{U}_{t_{i}}^{A}(s\ell(2,\mathbb{C}))}

\newcommand{\univg}{\mathcal{U}(\mathfrak{g})}
\newcommand{\uqg}{\mathcal{U}_{q}(\mathfrak{g})}
\newcommand{\utg}{\mathcal{U}_{t}(\mathfrak{g})}
\newcommand{\utga}{\mathcal{U}_{t}^{A}(\mathfrak{g})}

\newcommand{\utgg}{\mathcal{U}_{t}(\mathfrak{g}\oplus\mathfrak{g})}

\newcommand{\hlie}{\mathfrak{h}}

\newcommand{\otsltc}{\mathcal{O}_{t}(SL(2,\mathbb{C}))}

\newcommand{\otsltcazero}{\mathcal{O}_{t}^{A_{0}}(SL(2,\mathbb{C}))}

\newcommand{\otisltc}{\mathcal{O}_{t_i}(SL(2,\mathbb{C}))}

\newcommand{\otisltcazero}{\mathcal{O}_{t_i}^{A_{0}}(SL(2,\mathbb{C}))}
\newcommand{\tlotisltcazero}{\widetilde{\mathcal{O}^{A_{0}}_{t_i}(SL(2,\mathbb{C}))}}
\newcommand{\oqisltc}{\mathcal{O}_{q_i}(SL(2,\mathbb{C}))}

\newcommand{\otsutwo}{\mathcal{O}_{t}(SU(2))}

\newcommand{\otsutwoazero}{\mathcal{O}^{A_{0}}_{t}(SU(2))}
\newcommand{\oqsutwo}{\mathcal{O}_{q}(SU(2))}

\newcommand{\otisutwoa}{\mathcal{O}^{A}_{t_i}(SU(2))}
\newcommand{\otisutwoazero}{\mathcal{O}^{A_{0}}_{t_i}(SU(2))}
\newcommand{\oqisutwo}{\mathcal{O}_{q_i}(SU(2))}

\newcommand{\otslnc}{\mathcal{O}_{t}(SL(n+1,\mathbb{C}))}

\newcommand{\otslncazero}{\mathcal{O}^{A_{0}}_{t}(SL(n+1,\mathbb{C}))}

\newcommand{\oqsu}{\mathcal{O}_{q}(SU(n+1))}
\newcommand{\otsuazero}{\mathcal{O}^{A_{0}}_{t}(SU(n+1))}

\newcommand{\otg}{\mathcal{O}_{t}(G)}
\newcommand{\oqg}{\mathcal{O}_{q}(G)}
\newcommand{\otk}{\mathcal{O}_{t}(K)}
\newcommand{\oqk}{\mathcal{O}_{q}(K)}
\newcommand{\otga}{\mathcal{O}^{A}_{t}(G)}

\newcommand{\otka}{\mathcal{O}^{A}_{t}(K)}
\newcommand{\otgazero}{\mathcal{O}^{A_{0}}_{t}(G)}
\newcommand{\tlotgazero}{\widetilde{\mathcal{O}^{A_{0}}_{t}(G)}}
\newcommand{\otkazero}{\mathcal{O}^{A_{0}}_{t}(K)}

\newcommand{\otgnplus}{\mathcal{O}^{A_{0}}_{t}({G/N}^{+})}
\newcommand{\otgnminus}{\mathcal{O}^{A_{0}}_{t}({G/N}^{-})}

\newcommand{\conn}{\mathrm{Comp}}
\newcommand{\wt}{\mathrm{wt}}
\newcommand{\dom}{\mathrm{Dom}}

\newcommand{\lieg}{\mathfrak{g}}

\newtheorem{dfn}{Definition}[section]
\newtheorem{thm}[dfn]{Theorem}

\newtheorem{lmma}[dfn]{Lemma}

\newtheorem{ppsn}[dfn]{Proposition}

\newtheorem{crlre}[dfn]{Corollary}

\newtheorem{xmpl}[dfn]{Example}

\newtheorem{rmrk}[dfn]{Remark}

\newtheorem{xrcs}[dfn]{Exercise}

\newcommand{\bdfn}{\begin{dfn}\rm}
\newcommand{\bthm}{\begin{thm}}
\newcommand{\blmma}{\begin{lmma}}
\newcommand{\bppsn}{\begin{ppsn}}
\newcommand{\bcrlre}{\begin{crlre}}
\newcommand{\bxmpl}{\begin{xmpl}}
\newcommand{\brmrk}{\begin{rmrk}\rm}
\newcommand{\bxrcs}{\begin{xrcs}\rm\footnotesize}
\newcommand{\bprf}{\begin{proof}}

\newcommand{\edfn}{\end{dfn}}
\newcommand{\ethm}{\end{thm}}
\newcommand{\elmma}{\end{lmma}}
\newcommand{\eppsn}{\end{ppsn}}
\newcommand{\ecrlre}{\end{crlre}}
\newcommand{\exmpl}{\end{xmpl}}
\newcommand{\ermrk}{\end{rmrk}}
\newcommand{\exrcs}{\end{xrcs}}
\newcommand{\eprf}{\end{proof}}

\begin{document}
\author[1]{\sc Saikat Das\thanks{saikat20r@isid.ac.in, saikat1811@gmail.com}}
\author[1]{\sc Ayan Dey\thanks{ayan22r@isid.ac.in, studentayandey@gmail.com}}
\author[1]{\sc Arup Kumar Pal\thanks{arup@isid.ac.in, arupkpal@gmail.com}}
\affil[1]{Indian Statistical Institute, Delhi, INDIA}
\title{A triangular decomposition for the crystal lattice \\ of 
quantized function algebras}
\maketitle
 \begin{abstract}
We prove a triangular decomposition theorem for the lower crystal 
lattice $\mathcal{O}_{t}^{A_{0}}(G)$ of the quantized function 
algebra $\mathcal{O}_{t}^{}(G)$, where $G$ is a connected simply 
connected complex Lie group with Lie algebra $\mathfrak {g}$ of 
type $A_{n}$, $B_{n}$, $C_{n}$, $D_{n}$, $E_{6}$ or $E_{7}$. As 
a consequence, we prove the inclusion 
$\mathcal{O}_{t}^{A_{0}}(G)\subseteq\mathcal{O}_{t}^{A_{0}}(K)$ 
conjectured by Matassa \& Yuncken in these cases. 
Using this inclusion, we then prove that the notions of crystallized 
quantized function algebra given by Matassa \& Yuncken coincide with 
that of Giri \& Pal.
 \end{abstract}
{\bf AMS Subject Classification No.:}
17B37, 
22E46, 
20G42, 
46L67. 



 \section{Introduction}
The triangular decomposition 
$\lieg=\mathfrak{n}_{-}\oplus\hlie\oplus \mathfrak{n}_{+}$ of a  
Lie algebra $\lieg$ is a classical result that dates back to very 
early twentieth century in the works of Killing and Cartan, and later of Weyl
who put it in the form we know it today. The triangular decomposition 
of the Universal enveloping algebra $\univg$ 
\[
\univg=\mathcal{U}(\mathfrak{n}_{-})\otimes \mathcal{U}(\hlie)\otimes \mathcal{U}(\mathfrak{n}_{+})
\]
is a consequence of decomposition of $\lieg$ and the Poincare-Birkhoff-Witt (PBW) 
theorem. This result also goes back to the early twentieth century.
On the dual side, for the coordinate ring, one has a triangular decomposition 
or a factorization in terms of highest and lowest weight matrix coefficients.
This can be derived as a consequence of the algebraic Peter-Weyl theorem 
together with some basic representation theoretic results
and is also a classical result dating back to the first half of the 
twentieth century.
With the development of quantum groups in the mid to late 1980's,
their $q$-analogues, the quantized universal enveloping algebras and the
quantized functions algebras became primary objects of interest,
and many results on their classical counterparts found natural generalizations
to the quantized versions.
The triangular decomposition of $\univg$ was extended to the quantized universal 
enveloping algebras by Lusztig (\cite{Lus-1993aa}) and 
Rosso (\cite{Ros-1990aa}) during the early 1990's. 
The corresponding decomposition for the quantized coordinate ring $\otg$ can 
be traced back to Levendorski \& Soibelmen 
(see Proposition 3.2.1, \cite{LevSoi-1991aa} for the statement)
and 
Joseph (see Proposition 9.2.2, \cite{Jos-1995aa} for the statement with proof). 
After the development of the theory of crystal 
bases in the early 1990's, due independently to Kashiwara and Lusztig,
lattices in the coordinate ring over the subrings  
$A_{0}:=\{f(t)/g(t)\in \mathbb{Q}(t): g(0)\neq 0\}$
and 
$A:=\{t^{-m}p(t): m\in\mathbb{N}, p(t)\in \mathbb{Q}[t]\}$
of $\mathbb{Q}(t)$ became very important objects of study.
In 2009, Lusztig (Proposition~3.3, \cite{Lus-2009aa})
proved a triangular decomposition theorem for the lattice $\otga$ in $\otg$ over the Laurent polynomial ring $A$, and used his result to give a bound on the number of 
generators of the coordinate ring. 
Our main result in the present paper is a triangular decomposition theorem 
for the crystal lattice $\otgazero$ in  $\otg$ over the ring $A_{0}$, which
is the localization at zero of the polynomial subring 
$\mathbb{Q}[t]$ of $\mathbb{Q}(t)$.
We deal with deformations of Lie groups of types $A_{n}$, $B_{n}$, $C_{n}$, 
$D_{n}$, $E_{6}$ and $E_{7}$ in this paper. 
The remaining cases $E_{8}$, $F_{4}$ and $G_{2}$ have recently been
settled by the second author in a separate article (\cite{Dey-2026aa}).

The decomposition theorem (Theorem~\ref{th:Minuscule}) is proved in Section~3.2,
after developing the necessary material in Section~2 and Section~3.1.
In the remaining parts of the paper, we use the decomposition to
prove an inclusion 
$\mathcal{O}^{A_{0}}_{t}(G)\subseteq \mathcal{O}^{A_{0}}_{t}(K)$
between two different lattices conjectured
by Matassa \& Yuncken (\cite{MatYun-2023aa}), which then enables
us to establish a link between two different ways of defining 
the crystallization of a quantized function algebra given by
Giri \& Pal (\cite{GirPal-2022tv}) and by Matassa \& Yuncken (\cite{MatYun-2023aa}).
The former
used limits of scaled matrix entries of the fundamental representation
of $SU_{q}(n+1)$ to arrive at a set of generators and relations
that describes the crytallized algebra in the type $A_{n}$ case,
while the latter used Kashiwara's theory of crystal bases in a crucial way
to map a certain lattice inside $\otg$, for each $q\in(0,1)$, into the 
bounded operators on a Hilbert space and then defined the crystallized 
algebra as the algebra generated by the limits of these image
operators as $q\to 0$.
In a subsequent paper (\cite{GirPal-2024aa}), Giri \& Pal  gave a complete classification of all irreducible representations of the crystallized 
algebra $C(SU_{0}(n+1))$ defined by them. The consequences of the 
classification theorem indicated that their notion of crystallized
qunatized functions algebra should be closely related to the one defined by
Matassa \& Yuncken. In the last section of this paper, we use the 
inclusion
$\mathcal{O}^{A_{0}}_{t}(G)\subseteq \mathcal{O}^{A_{0}}_{t}(K)$
of lattices derived in Section~3 and the results in (\cite{GirPal-2024aa})
to establish equality of the two
notions of crystallization (Theorem~\ref{th:gpmy}).

For the necessary material on crystal basis theory, we refer the reader 
to the papers of Kashiwara (\cite{Kas-1990aa}, \cite{Kas-1991aa}, 
\cite{Kas-1993aa}) and the books by Joseph (Chapters 5--6, \cite{Jos-1995aa}),
Jantzen (Chapters 9--11, \cite {Jan-1996aa}) and Hong \& Kang
(Chapters 4-6, \cite{HonKan-2002ab}).
For a treatment on crystal lattices and crystal bases for the quantized 
function algebras, please refer to Kashiwara (\cite{Kas-1993aa}).
For various quantum function algebras and crystal lattices on
them, we will closely follow the notations of Matassa \& Yuncken
(\cite{MatYun-2023aa}) with one crucial difference that we elaborate 
in the next paragraph. So for example, $\otg$ will denote the quantized
function algebra for the group $G$ over the base field $\mathbb{Q}(t)$,
$\otgazero$ will denote the lower crystal lattice of the algebra $\otg$
with respect to the $\utg$-module structure on $\otg$ considered in
\cite{Igl-2006aa} and \cite{MatYun-2023aa}.
In Kashiwara's notation (\cite{Kas-1993aa}), $\otg$ and $\otgazero$ are
$A_{q}(\mathfrak{g})$ and $L(A_{q}(\mathfrak{g}))$ respectively, although
the $\utg$-module structure on $\otg$ there is different, resulting
in a slightly different lattice.
The space $\otga$ is the $A$-form  $A_{q}^{Q}(\mathfrak{g})$ in 
Kashiwara's notation, and $\otka$ is $\otga$ along with the $\ast$-structure
coming from the compact real form.

One important point on notations here.
In order to make a clear distinction between the cases when we work with
a deformation over the field of rational functions in an indeterminate
and when we work over the complex field, we will denote the 
the deformation parameter in the former case, which is the 
indeterminate, by `$t$' while in the latter case by `$q$'.
Also, in the former situation, we will work over the base 
field $\mathbb{Q}(t)$ rather than over $\mathbb{C}(t)$. 
Then the set up is identical to that in
\cite{Jos-1995aa}, \cite{Kas-1990aa} or \cite{Kas-1991aa}
and we can use their results freely.
In \cite{MatYun-2023aa}, Matassa \& Yuncken used the base field
$\mathbb{C}(t)$. The effect is that the spaces 
$\otg$, $\otk$ etc in their paper are complexifications
of the spaces denoted by the same symbols in this article. 
One also needs to work with the complexification of the symmetric
bilinear form (polarization) that we use here.
However, these do not have any effect on the final statement
on equality of the two crystallized algebras, as after passing
to $\oqk$ through the specialization maps, this difference disappears.

Let us next give a brief summary of the content. In the next section, we
fix the notations to be used for the rest of the article. Then
we describe a choice of basis that we make for each finite dimensional
irreducible highest weight module. It is important to fix this choice
in order to compare the two crystallizations, as we will see in the 
later sections. In Section~3, we first describe the crystal lattice and various
subspaces of $\otg$ used by Matassa \& Yuncken in \cite{MatYun-2023aa} in 
terms of these basis elements. We then prove a triangular decomposition 
theorem for the lower crystal lattice of the quantized function algebra
$\otg$ for all simple complex Lie algebras other than the exceptional classes
$E_{8}$, $F_{4}$ and $G_{2}$. We also give an alternative description of the space
$\mathcal{O}^{A_{0}}_t(G/N^{-})$. As a consequence of this and
the triangular decomposition, we derive the inclusion
$\mathcal{O}^{A_{0}}_{t}(G)\subseteq \mathcal{O}^{A_{0}}_{t}(K)$
for the above classes of Lie algebras.
In Section~4, we discuss the Matassa-Yuncken crystallization
in detail. In particular, we give a precise description of
the specialization (or evaluation) map in Subsection~4.2
and also give an alternate description of the crystallized algebra
where one applies the specialization map followed by the
Sobelman representation, avoiding the use of the restriction maps
which are subsumed in the Soibelmen representation at the
`$q$'-level. 
In the final section, we take up the type $A_{n}$ case.
It follows from the results in Section~3 that the lower crystal 
lattice $\otslncazero$ of the quantized algebra $\otslnc$ is the 
$A_{0}$-algebra generated by the matrix elements of the fundamental representation
$\varpi_{1}\equiv\one= (1,0,\ldots,0)$. In this section, we give a similar 
description of the $*$-algebra $\otsuazero$ in terms of the same
matrix elements after appropriately scaling them. This allows us 
to prove that the crystallized algebras defined in \cite{GirPal-2022tv} 
and \cite{MatYun-2023aa} are isomorphic.

 \section{Setup}

 \subsection{Preliminaries}
A symmetric bilinear form on a vector space $V$ will typically be denoted
by $(\cdot,\cdot)$, and the pairing between $V$ and its dual $V^{*}$ will
be denoted by $\langle\cdot,\cdot\rangle$. Given a finite dimensional
vector space $V$ and a nondegenerate symmetric bilinear form $(\cdot,\cdot)$ 
on $V$, let us denote by $v\mapsto v^{*}$
the map from $V$ to $V^{*}$ that gives an isomorphism between $V$ and $V^{*}$
through this form, i.e.\  such that 
$\langle v^{*},v^\prime\rangle=(v,v^\prime)$ for all $v,v^\prime\in V$.
Let $f\mapsto f_{*}$ denote the inverse of this isomorphism.

We will follow the standard notations from Lie theory. For example,
for a Lie algebra $\mathfrak{g}$ with a symmetrizable 
Cartan matrix $(\!(a_{i,j})\!)\in GL_{n}(\mathbb{Z})$, 
simple roots will be denoted by
$\alpha_{1},\ldots,\alpha_{n}$, with $\alpha^{\vee}_{1},\ldots,\alpha^{\vee}_{n}$ 
the corresponding coroots such that 
$\langle \alpha^{\vee}_{i},\alpha^{}_{j}\rangle=a_{i,j}$.
The fundamental weights will be denoted by 
$\varpi_{1},\ldots,\varpi_{n}$.
The set of highest weights will be denoted by $P_{+}$ and
for $\Lambda\in P_{+}$, $P(\Lambda)$ will denote the set of weights
in the irreducible highest weight module $V(\Lambda)$.
The sum of the fundamental weights will be denoted by $\rho$.
The quantities $d_{i}$ will denote $\frac{1}{2}(\alpha_{i},\alpha_{i})$
where $(\cdot,\cdot)$ is a symmetric nondegenerate bilinear form, normalized
so as to make each $d_{i}$ an integer and the $d_{i}$'s coprime. 
For an integer linear combination $\mu=\sum_{i}m_{i}\alpha_{i}$, 
we will denote by $K^{\mu}$ the product $K_{1}^{m_{1}}\cdots K_{n}^{m_{n}}$.

Let $\mathfrak{g}$ be a simple complex Lie algebra with 
Cartan matrix $(\!(a_{i,j})\!)$ and $G$ be the unique
connected simply connected Lie group with $\mathfrak{g}=\text{Lie}\, G$.
The quantized universal enveloping algebra $\utg$
is an algebra over $\mathbb{Q}(t)$ generated by $K_i$, $K_i^{-1}$, $E_{i}$ 
and $F_{i}$ for $i=1,\ldots,n$ obeying the following relations:
\begin{IEEEeqnarray}{rClClrCrCl}
   K_i K_i^{-1} &=&K_i^{-1} K_i &=& 1,&&\qquad \qquad&
   K_i K_j &=& K_j K_i  , \label{eq:uqrel-1}\\ 
   K_i E_{j} K_i^{-1} &=& \IEEEeqnarraymulticol{3}{l}{t^{d_{i}a_{i j}} E_{j} }&&&
   K_i F_{j} K_i^{-1} &=& t^{-d_{i}a_{i j}} F_{j}, \label{eq:uqrel-2}\\
   E_{i} F_{j}-F_{j} E_{i} &=& \IEEEeqnarraymulticol{8}{l}{\delta_{i j} \frac{K_i-K_i^{-1}}{t^{d_{i}}-t^{-d_{i}}}},\label{eq:uqrel-3}\\
   \IEEEeqnarraymulticol{8}{r}{
   \sum_{k=0}^{1-a_{i j}}(-1)^k\left[\begin{array}{c}
1-a_{i j} \\
k
\end{array}\right]_{t^{d_{i}}} 
E_{i}^{1-a_{i j}-k} E_{j} E_{i}^k 
}
&=& 0  \quad (i \neq j) \label{eq:uqrel-4}\\
\IEEEeqnarraymulticol{8}{r}{
\sum_{k=0}^{1-a_{i j}}(-1)^k\left[
\begin{array}{c}
1-a_{i j} \\
k
\end{array}\right]_{t^{d_{i}}} F_{i}^{1-a_{i j}-k} F_{j} F_{i}^k }
&=& 0 \quad (i \neq j),\label{eq:uqrel-5}
\end{IEEEeqnarray}
where $\text{diag}(s_{1},\ldots,d_{n})$ is the symmetrizer
of the Cartan matrix $(\!(a_{i,j})\!)$.
Here the $s$-binomial coefficient 
${\genfrac{[}{]}{0pt}{}{m}{k}}_{s}$
denotes $\frac{(m)_{s}!}{(k)_{s}!(m-k)_{s}!}$,
with the $s$-numbers $(k)_{s}$ being given by $\frac{s^{k}-s^{-k}}{s-s^{-1}}$.
The Hopf structure is given by the following maps:
\begin{IEEEeqnarray}{rClrClrCl}
\Delta(K_i^{}) &=& K_{i}^{}\otimes K_{i}^{},\quad&
    \Delta(E_{i}^{}) &=& E_{i}^{}\otimes K_{i}^{-1}+ 1\otimes E_{i}^{},\quad&
    \Delta(F_{i}^{}) &=& F_{i}^{}\otimes 1 + K_{i}^{}\otimes F_{i}^{},
        \label{eq:uqrel-6}\\
S(K_{i}^{}) &=& K_{i}^{-1}, &
    S(E_{i}^{}) &=& -E_{i}^{}K_{i}^{},&
       S(F_{i}^{}) &=& - K_{i}^{-1}F_{i}^{},\label{eq:uqrel-7}
\end{IEEEeqnarray}
and the compact real form on $\utg$ is the following $*$-structure on it:
\begin{IEEEeqnarray}{rClrClrCl}
\label{eq:cptreal}
K_{i}^{*} &=& K_{i}^{}, \qquad &
E_{i}^{*} &=& t^{d_{i}}F_{i}^{}K_{i}^{-1},\qquad &
F_{i}^{*} &=& t^{-d_{i}}K_{i}^{}E_{i}^{}.
\end{IEEEeqnarray}

For $q\in \mathbb{C}$, $q\neq 0,1$, the quantized universal 
enveloping algebra $\uqg$ is a Hopf algebra over $\mathbb{C}$ 
generated by the same relations (\ref{eq:uqrel-1}--\ref{eq:uqrel-7}) 
as above with the indeterminate $t$ replaced by the complex number $q$. 
The quantized function algebra $\oqg$ is the subspace of $\uqg^{*}$ 
spanned by the matrix elements of the finite dimensional irreducible
highest weight modules of $\uqg$. 
When $q\in (0,\infty)$, one has the compact real form $\oqk$, which
is $\oqg$ together with the involution 
\begin{IEEEeqnarray}{rClrClrCl}
\label{eq:cptreal-q}
K_{i}^{*} &=& K_{i}^{}, \qquad &
E_{i}^{*} &=& q^{d_{i}}F_{i}^{}K_{i}^{-1},\qquad &
F_{i}^{*} &=& q^{-d_{i}}K_{i}^{}E_{i}^{}.
\end{IEEEeqnarray}
We will sometimes use $t_{i}$ and $q_{i}$ to denote $t^{d_{i}}$ 
and $q^{d_{i}}$ respectively.
Let $n$ be the rank of $\lieg$ and let $k$ be the length of the longest
element  of the Weyl group. Then the Hilbert space 
$\ell^{2}(\mathbb{N})^{\otimes k}\otimes \ell^{2}(\mathbb{Z})^{\otimes n}$
will be denoted by $\mathcal{H}_{Soi}$ and the Soibelmen representation of
$C(K_{q})$ on this Hilbert space by $\psi^{(q)}_{Soi}$.
For any Hilbert space $\mathcal{H}$, the space of bounded linear operators
on $\mathcal{H}$ will be denoted by $\mathcal{L}(\mathcal{H})$.

\brmrk
The algebra and the coalgebra structure and the compact
real form given above are identical to those used 
in \cite{Kas-1991aa} and \cite{MatYun-2023aa}. The coproduct
map $\Delta$ and the antipode map $S$ here, in particular, 
are same as the coproduct
$\Delta_{-}$ and the corresponding antipode map used by 
Kashiwara while treating lower crystal
basis. These maps are different in different
articles or books and can be a big source of confusion, though
it is usually possible to pass from one to the other by using 
a Hopf algebra automorphism, in other words, by passing on to a different
set of generators.
\ermrk

\subsection{Choice of basis}
In order to deal with the two notions of crystallization given by
Giri \& Pal and Matassa \& Yuncken, it will be convenient to make 
a choice of basis for each finite dimensional irreducible highest weight module
of $\utg$. We will do that as a certain lifting of the crystal 
base for each such module.

To begin with, let us  choose a highest weight vector 
$u^{\Lambda}_{\Lambda}$ in the highest weight module $V(\Lambda)$ for each 
$\Lambda\in P_{+}$. Let $L_{u}(\Lambda)$ be the corresponding
lower crystal lattice, $(\cdot,\cdot)_{u}$ the polarization on $V(\Lambda)$,
and $\{u^{\Lambda}_{i}:1\leq i\leq\dim\Lambda\}$ the lower global basis
with $u^{\Lambda}_{1}=u^{\Lambda}_{\Lambda}$.

Let $\omega_{0}$ be the longest element of the Weyl group $W$.
Let us denote 
$(u^{\Lambda}_{\omega_{0}\Lambda}, u^{\Lambda}_{\omega_{0}\Lambda})$ 
by $c_{\Lambda}$.
Observe that $c_{\Lambda}\in 1+tA_{0}$. In particular, it is
invertible in $A_{0}$.
Define a bilinear form on $V(\Lambda)^{*}\cong V(-\omega_{0}\Lambda)$
by
\begin{equation}
\label{eq:dual-pol}
((z_{1})^{*}, (z_{2})^{*})_{\tilde{u}}:= c_{\Lambda}^{-1} 
     t^{(\Lambda, 2\rho)}(K^{-2\rho}z_{1}, z_{2})_{u},
          \quad z_{1},z_{2}\in V(\Lambda).
\end{equation}
Note that this is same as the form defined on $V(\Lambda)^{*}$
in \cite{StoDij-1999aa} except for the normalizing factor
$t^{(\Lambda, 2\rho)}$ that is needed in our case to ensure 
that certain vectors we consider belong to the crystal lattices
of the respective spaces.
It is easy to see that the form defined above is symmetric, 
nondegenerate and makes the module $V(-\omega_{0}\Lambda)$ unitarizable. 
The vector $(u^{\Lambda}_{\omega_{0}\Lambda})^{*}$ is of highest weight 
in $V(-\omega_{0}\Lambda)$. 
Let us take the scaled vector
\begin{equation} \label{eq:v_1}
{\tilde{u}}^{-\omega_{0}\Lambda}_{1}:=
t^{(\omega_{0}\Lambda-\Lambda,\rho)}(u^{\Lambda}_{\omega_{0}\Lambda})^{*}
\end{equation}
as our choice of highest weight vector and let 
$L_{\tilde{u}}(-\omega_{0}\Lambda)$ be the corresponding lower crystal 
lattice. Denote by 
$\{{\tilde{u}}^{-\omega_{0}\Lambda}_{i}: 1\leq i\leq \dim \Lambda\}$
the lower global basis for $V(-\omega_{0}\Lambda)$ corresponding to this choice.
Note that by Theorem~6.2.2, \cite{HonKan-2002ab}, this is the unique lifting
of the corresponding crystal basis that is invariant under the bar map `--'.
Observe also that
$({\tilde{u}}^{-\omega_{0}\Lambda}_{1}, 
         {\tilde{u}}^{-\omega_{0}\Lambda}_{1})_{\tilde{u}}=1$, 
which implies that the form $(\cdot,\cdot)_{\tilde{u}}$ is the polarization on 
$V(-\omega_{0}\Lambda)$.

If $\Lambda=-\omega_{0}\Lambda$, there exists a $\utg$-module
isomorphism
\begin{equation}\label{eq:iso-to-dual}
T: V(\Lambda)\to V(-\omega_{0}\Lambda)\cong V(\Lambda)^{*}
\end{equation}
such that $T(u^{\Lambda}_{1})={\tilde{u}}^{-\omega_{0}\Lambda}_{1}$,
$(T(z_{1}), T(z_{2}))_{\tilde{u}}=(z_{1}, z_{2})_{u}$ and
$T(L_{u}(\Lambda))=L_{\tilde{u}}(-\omega_{0}\Lambda)$.

\bppsn
\label{ppsn:dual-lattice}
For each $1\leq k\leq m:=\dim\Lambda$, let
\[
w^{-\omega_{0}\Lambda}_{k}=
  t^{(\wt(u^{\Lambda}_{m+1-k}) - \Lambda, \rho)} 
  (u^{\Lambda}_{m+1-k})^{*}\in V(\Lambda)^{*}\cong V(-\omega_{0}\Lambda).
\]
  Then one has 
$w^{-\omega_{0}\Lambda}_{k}\in L_{\tilde{u}}(-\omega_{0}\Lambda)$
and
\[
L_{\tilde{u}}(-\omega_{0}\Lambda)=
\text{$A_{0}$-span of }\left\{w^{-\omega_{0}\Lambda}_{k}: 1\leq k \leq m \right\}.
\]
Furthermore, for all $\mu\in P(-\omega_{0}\Lambda)$, one has
\begin{equation}
\label{eq:wt-sp}
L_{\tilde{u}}(-\omega_{0}\Lambda)_{\mu} = \mathrm{Span}_{A_0} 
 \left\{w^{-\omega_{0}\Lambda}_{k} : 
    \wt\left(v^{\Lambda}_{m-k+1}\right) = -\mu\right\}.
\end{equation}
\eppsn
\bprf
For the course of proof of this result, we will denote  by $L$
the $A_{0}$-span of 
\[
\left\{w^{-\omega_{0}\Lambda}_{k}: 1\leq k \leq m \right\}.
\]
Let $\mu$ be the weight of $u^{\Lambda}_{k}$, where $1\leq k\leq \dim\Lambda$. 
Then
\begin{IEEEeqnarray*}{rCl}
 ((u^{\Lambda}_{k})^{*},(u^{\Lambda}_{k})^{*})
  &=&
  c_{\Lambda}^{-1} 
     t^{(\Lambda, 2\rho)}(K^{-2\rho}u^{\Lambda}_{k}, u^{\Lambda}_{k})\\
  &=&
  c_{\Lambda}^{-1} 
     t^{(\Lambda-\mu, 2\rho)}(u^{\Lambda}_{k}, u^{\Lambda}_{k}).
\end{IEEEeqnarray*}
By Proposition~6.1.2, \cite{Jos-1995aa}, 
it follows  that 
$w^{-\omega_{0}\Lambda}_{k}\in L(-\omega_{0}\Lambda)$ for $1\leq k\leq m$.
Therefore  $L\subseteq L_{\tilde{u}}(-\omega_{0}\Lambda)$.
In particular, we also have 
$L + t L_{\tilde{u}}(-\omega_{0} \Lambda) \subseteq 
     L_{\tilde{u}}(-\omega_{0} \Lambda)$.

Next, observe that the set 
$\left\{w^{-\omega_{0}\Lambda}_{k} + t L_{\tilde{u}}(-\omega_{0} \Lambda) : 
 1\leq k \leq m\right\}$ 
 is orthonormal in the $\mathbb{Q}$-vector space 
 $L_{\tilde{u}}(-\omega_{0}\Lambda)/t L_{\tilde{u}}(-\omega_{0}\Lambda)$
 which has dimension $m$. Therefore
 $L_{\tilde{u}}(-\omega_{0}\Lambda)\subseteq L+tL_{\tilde{u}}(-\omega_{0}\Lambda)$.
The equality $L=L_{\tilde{u}}(-\omega_{0}\Lambda)$ now follows from 
Nakayama's Lemma (see for example Proposition~2.6, 
\cite{AtiMac-1969aa}) over $A_{0}$. 

 A similar argument gives us the equality (\ref{eq:wt-sp}).
\eprf

\brmrk
\label{rm:glb-dual}
Observe that for each $k$, 
$w^{-\omega_{0}\Lambda}_{k}$ is $(u^{\Lambda}_{m-k+1})^{*}$
normalized so that 
$(w^{-\omega_{0}\Lambda}_{k}, w^{-\omega_{0}\Lambda}_{k})$ is
in $1+tA_{0}$. In particular, 
$(w^{-\omega_{0}\Lambda}_{k}, w^{-\omega_{0}\Lambda}_{k})$
is a unit (i.e.\ invertible) in $A_{0}$.
Therefore if $u^{\Lambda}_{m-k+1}$ comes from a one-dimensional weight space,
then $w^{-\omega_{0}\Lambda}_{k}=c. {\tilde{u}}^{-\omega_{0}\Lambda}_{k}$
for some unit $c$ with $c^2\in 1+tA_{0}$.
In particular, $w^{-\omega_{0}\Lambda}_{m}=c.{\tilde{u}}^{-\omega_{0}\Lambda}_{m}$
for some $c\in \pm 1+tA_{0}$.
\ermrk
\bcrlre
Assume $\Lambda\neq 0$. Then
\[
\text{$A_{0}$-span of }\left\{(u^{\Lambda}_{k})^{*}: 1\leq k \leq m\right\}
 \subsetneqq L_{\tilde{u}}(-\omega_{0}\Lambda).
\]
\ecrlre
\bprf
Since $\Lambda-\wt(u^{\Lambda}_{k})$ is a linear combination of
$\alpha_{i}$'s with nonnegative integer coefficients and
$\rho$ is the sum of all fundamental weights, it follows that
$(\Lambda-\wt(u^{\Lambda}_{k}),\rho)$ is a nonnegative integer.
Therefore it follows from the above proposition that
$(u^{\Lambda}_{k})^{*}\in L_{\tilde{u}}(-\omega_{0}\Lambda)$.

Assume $L_{\tilde{u}}(-\omega_{0}\Lambda)=\text{$A_{0}$-span of }
        \left\{(u^{\Lambda}_{k})^{*}: 1\leq k \leq m\right\}$.
Since $\Lambda\neq 0$, we have $\text{dim}\,V(\Lambda)>1$
and hence $V(\Lambda)$ has a weight vector with weight 
$\mu\neq\Lambda$. Then we have, for some $k$,
\[
t^{(\mu-\Lambda,\rho)}(u^{\Lambda}_{k})^{*} 
   = \sum_{j=1}^{m}f_{j}(t)(u^{\Lambda}_{j})^{*},
\]
where $f_{j}(t)\in A_{0}$ for all $j$.
Taking the pairing of both sides with $v^{\Lambda}_{k}$
and noting that $(u^{\Lambda}_{i},u^{\Lambda}_{j})\in \delta_{i,j}+tA_{0}$,
one arrives at a contradiction.
\eprf

\paragraph{Choice of basis.}\text{}\\ 
We now make the following choice of global basis for each
$\Lambda\in P_{+}$.
Let $\Lambda=\sum_{i}m_{i}\varpi_{i}$  and let 
   $-\omega_{0}\Lambda=\sum_{i}n_{i}\varpi_{i}$.
If $\Lambda \geq -\omega_{0}\Lambda$ in lexicographic ordering
(identifying $\Lambda$ and $-\omega_{0}\Lambda$ with the tuples
$(m_{i})$ and $(n_{i})$ respectively), then we will take
the crystal lattice to be $L(\Lambda):=L_{u}(\Lambda)$,
the global basis $v^{\Lambda}_{i}:=u^{\Lambda}_{i}$
and polarization $(\cdot,\cdot):=(\cdot,\cdot)_{u}$.
If $\Lambda < -\omega_{0}\Lambda$, then we will work
with the crystal lattice $L(\Lambda):=L_{\tilde{u}}(\Lambda)$,
the global basis $v^{\Lambda}_{i}:={\tilde{u}}^{\Lambda}_{i}$
and polarization $(\cdot,\cdot):=(\cdot,\cdot)_{\tilde{u}}$.

\section{The crystal lattice}

\subsection{Subspaces of the coordinate ring}
For $f\in V(\Lambda)^{*}$ and $v\in V(\Lambda)$, we will denote by
$C^{\Lambda}_{f,v}$ the matrix element given by 
$C^{\Lambda}_{f,v}(a)=\langle f, av\rangle$.
We will work with a fixed polarization on $V(\Lambda)$
and as mentioned earlier, $w^{*}$ will denote the functional $v\mapsto (w,v)$
so that $C^{\Lambda}_{w^{*},v}$ will denote 
the matrix element $a\mapsto \langle w^{*},v\rangle =(w,av)$.
Occasionally we will also use the notation $C^{\Lambda}_{w,v}$
to denote the same matrix element.
When we  work with a basis $\{v^{\Lambda}_{i}:1\leq i\leq \dim\Lambda\}$,
we will usually write $C^{\Lambda}_{i,j}$
in place of $C^{\Lambda}_{(v^{\Lambda}_{i})^{*},v^{\Lambda}_{j}}
 \equiv C^{\Lambda}_{v^{\Lambda}_{i},v^{\Lambda}_{j}}$.

Let us start by recalling the set up and the crystal lattice on $\otg$
described in Iglesias \cite{Igl-2006aa} that was used by Matassa \& Yuncken
in \cite{MatYun-2023aa}. The algebra $\utg$ has a two commuting actions 
$\mathcal{L}$ and  $\mathcal{R}$ on $\utg^{*}$ given by
\[
(\mathcal{L}(a)f)(b):= f(a^{*}b), \qquad (\mathcal{R}(a)f)(b):= f(ba).
\]
These translate to the following actions on the matrix 
elements $C^{\Lambda}_{v,w}$, where 
$\Lambda\in P_{+}$ and $v,w\in V(\Lambda)$:
\[
\mathcal{L}(a)C^{\Lambda}_{v,w}= C^{\Lambda}_{av,w},\qquad
  \mathcal{R}(a)C^{\Lambda}_{v,w}= C^{\Lambda}_{v,aw},
\]
which give $\otg$ a $\utg\otimes\utg\cong \utgg$-module
structure given by
\[
((a\otimes b)\cdot f)(c):=f(a^{*}cb) 
   = (\mathcal{L}(a)\mathcal{R}(b)f)(c).
\]
One then has the the direct sum decomposition as a $\utgg$-module
\[
\otg=\bigoplus_{\Lambda\in P_{+}} \mathbb{Q}(t)\text{-span of }
   \left\{C^{\Lambda}_{v^{\Lambda}_{i},v^{\Lambda}_{j}}: 
        1\leq i,j\leq\dim\Lambda\right\}
\]
and the resulting identification
\begin{equation}\label{eq:otg-identification1}
\otg\cong \bigoplus_{\Lambda\in P_{+}}
   \left(V(\Lambda)\otimes V(\Lambda)\right)
\end{equation}
through the map 
\begin{equation}\label{eq:otg-identification2}
C^{\Lambda}_{v,w}\mapsto v\otimes w.
\end{equation}
Consequently, we have
\begin{equation}\label{eq:otg-identification3}
\otgazero\cong \bigoplus_{\Lambda\in P_{+}}
   \left(L(\Lambda)\otimes_{A_{0}} L(\Lambda)\right).
\end{equation}
By the choices of basis we have made in the previous section, we have
\[
L(\Lambda)=A_{0}\text{-span of }\{v^{\Lambda}_{i}: 
       \Lambda\in P_{+}, 1\leq i\leq \dim\Lambda\}.
\]
Thus via the map (\ref{eq:otg-identification2}), we have
\begin{IEEEeqnarray}{rCl}
 \otgazero
 &=& A_{0}\text{-span of }\Bigl(
     \left\{ C^{\Lambda}_{(v^{\Lambda}_{i})^{*}, v^{\Lambda}_{j}}: 
       \Lambda\in P_{+},  1\leq i,j\leq \dim\Lambda\right\} .
   \label{eq:otgazero}
\end{IEEEeqnarray}

\brmrk \label{rm:tlotgazero1}
Let $\{(v^{\Lambda}_{i})^{o}: 1\leq i\leq \dim\Lambda\}$ be the basis of 
$V(\Lambda)^{*}\cong V(-\omega_{0}\Lambda)$ dual to the basis
$\{v^{\Lambda}_{i}: 1\leq i\leq \text{dim\,}\Lambda\}$ of $V(\Lambda)$.
Recall  (Definition~3.2, \cite{MatYun-2023aa}) that 
$\otgnplus$ is the $A_{0}$-subalgebra of $\otg$ generated by
$\{C^{\Lambda}_{(v^{\Lambda}_{i})^{o},v^{\Lambda}_{1}}: 
     \Lambda\in P_{+}, 1\leq i\leq \dim\Lambda\}$,
the space $\otgnminus$ is the $A_{0}$-subalgebra of $\otg$ generated by
$\{S(C^{\Lambda}_{(v^{\Lambda}_{1})^{o},v^{\Lambda}_{i}}): 
     \Lambda\in P_{+}, 1\leq i\leq \dim\Lambda\}$,
     and $\otkazero$ is the $A_{0}$-subalgebra of $\otg$
     generated by $\otgnplus$ and $\otgnminus$.
 Note that $\otgnplus\subseteq \otgazero$ and 
 $\otgnminus\subseteq (\otgazero)^{*}$. Therefore 
 $\otkazero$ is contained in the $\ast$-algebra $\tlotgazero$ generated
 by $\otgazero$.
\ermrk

Note that for each highest weight $\Lambda$, Matassa \& Yuncken 
(\cite{MatYun-2023aa}) chose and fixed a lifting of the 
crystal basis. Here we have used a lower global basis 
$\{v^{\Lambda}_{i}:i\}$. This will allow us to use
the additional properties enjoyed by a global basis.

Observe that
\begin{IEEEeqnarray*}{rCl}
  S(C^{\Lambda}_{(v^{\Lambda}_{1})^{o},v^{\Lambda}_{i}})(a)
  &=&
  C^{\Lambda}_{(v^{\Lambda}_{1})^{o},v^{\Lambda}_{i}}(Sa)\\
  &=&
  \langle (v^{\Lambda}_{1})^{o}, (Sa) v^{\Lambda}_{i}\rangle\\
  &=&
  \langle a(v^{\Lambda}_{1})^{o},  v^{\Lambda}_{i}\rangle\\
  &=&
  C^{-\omega_{0}\Lambda}_{v^{\Lambda}_{i}, (v^{\Lambda}_{1})^{o}}(a).
\end{IEEEeqnarray*}
Therefore  the $A_{0}$-algebra $\otgnminus$ is
generated by
$\{C^{-\omega_{0}\Lambda}_{v^{\Lambda}_{i}, (v^{\Lambda}_{1})^{o}}: 
     \Lambda\in P_{+}, 1\leq i\leq \dim\Lambda\}$.

\bppsn
\label{ppsn:rplus}
The algebra $\mathcal{O}^{A_{0}}_t(G/N^{+})$ is the $A_{0}$-subalgebra of $\otg$
generated by the elements
$\{C^{\Lambda}_{i,1}:  \Lambda\in P_{+}, 1\leq i\leq \dim\,\Lambda\}$.
\eppsn
\bprf
Take a $\Lambda\in P_{+}$. Let $m$ be the dimension of $V(\Lambda)$.
Let $\{(v^{\Lambda}_{i})^{o}: 1\leq i\leq m\}$ be dual basis of 
$V(\Lambda)^{*}\cong V(-\omega_{0}\Lambda)$ as before.
Then there is a matrix $(\!(p^{\Lambda}_{i,j}(t))\!)$ 
in $GL_{m}(A_{0})$ with 
$p_{i,j}(t)\in\delta_{i,j}+tA_{0}$
such that
\begin{IEEEeqnarray}{rCl}
(v^{\Lambda}_{i})^{*}
 &=&
 \sum_{j=1}^{m}p^{\Lambda}_{i,j}(t) (v^{\Lambda}_{j})^{o}, \quad 1\leq i\leq m.
 \label{eq:v^o=v^*}
\end{IEEEeqnarray}
Note that in such a case, the inverse of the matrix 
$(\!(p^{\Lambda}_{i,j}(t))\!)$ is again 
of the same form as $(\!(p^{\Lambda}_{i,j}(t))\!)$.
Since $\mathcal{O}^{A_{0}}_t(G/N^{+})$ is 
the $A_{0}$-subalgebra of $\otg$ generated by
$\{C^{\Lambda}_{(v^{\Lambda}_{i})^{o},v^{\Lambda}_{1}}: 
     \Lambda\in P_{+}, 1\leq i\leq m\}$,
the result follows by equation~(\ref{eq:v^o=v^*}).
\eprf

We will denote by  $R^{A_{0}}_{+}$ the $A_{0}$-span of 
$\{C^{\Lambda}_{v^{\Lambda}_{i},v^{\Lambda}_{1}}: 
      \Lambda\in P_{+}, 1\leq i\leq \dim\Lambda\}$.
and by $R^{A_{0}}_{-}$ the $A_{0}$-span of 
$\{C^{\Lambda}_{v^{\Lambda}_{i},v^{\Lambda}_{\omega_{0}\Lambda}}: 
    \Lambda\in P_{+}, 1\leq i\leq \dim\Lambda\}$.
These are the crystal analogues of the spaces $R_{q}(G/N^{+})$ and
$R_{q}(G/N^{-})$ defined in \cite{Jos-1995aa}. It is known
(see 9.1.6, page 264, \cite{Jos-1995aa}) that $R_{q}(G/N^{+})$ and
$R_{q}(G/N^{-})$ are $\mathbb{Q}(t)$-algebras. Let us next prove that
an analogous statement holds for $R^{A_{0}}_{+}$ and $R^{A_{0}}_{-}$.

\bppsn
\label{ppsn:rplusminus-alg}
$R^{A_{0}}_{+}$ and $R^{A_{0}}_{-}$ are algebras over $A_{0}$.
\eppsn
\bprf
Take two matrix elements $C^{\Lambda}_{r,1}$ and $C^{\Omega}_{s,1}$ 
in $R^{A_{0}}_{+}$. We will show that their product belongs to $R^{A_{0}}_{+}$.
The simple module $V(\Lambda+\Omega)$ appears in the direct sum decomposition
of $V:=V(\Lambda)\otimes V(\Omega)$ into irreducible modules with multiplicity one,
called the Cartan component of $V(\Lambda)\otimes V(\Omega)$,
and one has
\[
V(\Lambda)\otimes V(\Omega)\cong 
 V(\Lambda+\Omega)\bigoplus \left(\bigoplus_{\substack{\Gamma\in P_{+}(V)\\ \Gamma < \Lambda+\Omega}} V(\Gamma)\right)
\]
Thus the weight space of weight $\Lambda+\Omega$ is one dimensional, 
generated by the vector $v^{\Lambda}_{1}\otimes v^{\Omega}_{1}$. 
Therefore there is a $\mathcal{U}_{t}(g)$-module isomorphism 
$T_{0}:$ $V(\Lambda+\Omega)\bigoplus 
\left(\bigoplus_{\substack{\Gamma\in P_{+}(V)\\ \Gamma < \Lambda+\Omega}} 
V(\Gamma)\right)\to V$ 
whose restriction is a $A_{0}$-linear isomorphism between 
$L(\Lambda+\Omega)\bigoplus(\bigoplus_
{\Gamma\neq \Lambda+\Omega} L(\Gamma))$ 
and $L(\Lambda)\otimes L(\Omega)$. Take the highest weight vector 
$v^{\Lambda+\Omega}_{1}$  in $V(\Lambda+\Omega)$ chosen earlier.
Then it follows that
$T_{0}(v^{\Lambda+\Omega}_{1}\oplus 0)= c.v^{\Lambda}_{1}\otimes v^{\Omega}_{1}$
where $c$ is a unit in $A_{0}$. 
Hence, by renormalizing $T_{0}$ and composing it with the embedding 
$V(\Lambda+\Omega)\to V(\Lambda+\Omega)\bigoplus 
(\bigoplus_{\Gamma\neq \Lambda+\Omega} V(\Gamma))$, 
we conclude that there exists a $\mathcal{U}_{t}(g)$-module 
morphism $T: V(\Lambda+\Omega)\to V$ such that 
$T(v^{\Lambda+\Omega}_{1})= v^{\Lambda}_{1}\otimes v^{\Omega}_{1}$ and 
$T(L(\Lambda+\Omega))\subseteq L(\Lambda)\otimes L(\Omega)$.  
 For any $a\in \mathcal{U}_{t}(g)$, we then have 
\begin{IEEEeqnarray*}{rCl}
  C^{\Lambda}_{v_{r}^{\Lambda},v_{1}^{\Lambda}}.
      C^{\Omega}_{v_{s}^{\Omega},v_{1}^{\Omega}}(a) &=& 
    C^{V}_{v_{r}^{\Lambda}\otimes v_{s}^{\Omega}, 
          v_{1}^{\Lambda}\otimes v_{1}^{\Omega}}(a)\\
   &=&
   \langle (v^{\Lambda}_{r}\otimes v^{\Omega}_{s})^{*}, 
        a (v^{\Lambda}_{1}\otimes v^{\Omega}_{1})\rangle\\
   &=&
   \langle (v^{\Lambda}_{r}\otimes v^{\Omega}_{s})^{*}, 
        a T(v^{\Lambda+\Omega}_{1})\rangle\\
   &=&
   \langle (v^{\Lambda}_{r}\otimes v^{\Omega}_{s})^{*}, 
        T(av^{\Lambda+\Omega}_{1})\rangle.
\end{IEEEeqnarray*}
Now observe that 
$T^{tr}(v^{\Lambda}_{r}\otimes v^{\Omega}_{s})^{*}\in 
Hom_{A_{0}}( L(\Lambda+\Omega), A_{0})$, so that
  there exist $c_{k}\in A_{0}$ such that 
  $T^{tr}(v^{\Lambda}_{r}\otimes v^{\Omega}_{s})^{*}
   = \sum_{k} c_{k}. (v^{\Lambda+\Omega}_{k})^{*}$. 
Therefore, we must have 
\[
C^{\Lambda}_{r,1}.C^{\Omega}_{s,1}= 
\sum_{k} c_{k}. C^{\Lambda+\Omega}_{v_{k}^{\Lambda+\Omega},v^{\Lambda+\Omega}_{1}}.
\] 
Thus $R^{A_{0}}_{+}$ is an $A_{0}$-algebra.

For $R^{A_{0}}_{-}$, use a crystal morphism $T$ that sends 
$v^{\Lambda+\Omega}_{\omega_{0}. (\Lambda+\Omega)}$
to $v^{\Lambda}_{\omega_{0}\Lambda}\otimes v^{\Omega}_{\omega_{0}\Omega}$. 
The rest of the argument is identical.
\eprf
Thus we have $\mathcal{O}^{A_{0}}_t(G/N^{+})=R^{A_{0}}_{+}$.
We will see in the next subsection that the algebra $\mathcal{O}^{A_{0}}_t(G/N^{-})$
is strictly larger than  $R^{A_{0}}_{-}$.

 \subsection{Decomposition of the crystal lattice}
In this subsection, we will prove a triangular decomposition 
for the crystal lattice $\otgazero$. For the decomposition
of the coordinate ring, we refer the reader to Joseph
(Proposition 9.2.2, \cite{Jos-1995aa}), and for an $A$-lattice
version of the decomposition, to Lusztig (Proposition~3.3, \cite{Lus-2009aa}).
 We deal with the case of all simple complex 
Lie algebras except the three exceptional classes $E_{8}$, $F_{4}$ and $G_{2}$.
This decomposition helps understand the structure of the spaces $\otgazero$ and 
$\otkazero$ better and we will exploit that to give a proof of the inclusion
$\otgazero\subseteq \otkazero$ in the next subsection.

Given a weight $\mu\in P$, we denote by $\dom(\mu)$ the unique dominant weight 
in the Weyl orbit of $\mu$. Let us start with two lemmas. 
\blmma \label{lm:PRV}
Let $\lieg$ be a simple complex Lie algebra and 
let $\Omega\neq 0$ be a dominant weight. Take a weight of the form 
$\mu:=\omega.\Omega$, where $\omega$ is an element of the Weyl group $W$. 
Then there are dominant weights $\Lambda, \Gamma$   such that $V(\Omega)$ 
appears as a direct summand of $V:=V(\Lambda)\otimes V(-\omega_{0}.\Gamma)$ 
with multiplicity 1 and $\Lambda-\Gamma= \mu$. 
\elmma
\begin{proof}
    Choose $\Gamma\in P_{+}$ such that $\mu+\Gamma \in P_{+}$. 
    Take $\Lambda:= \mu+\Gamma$, and the rest follows from 
    the PRV Theorem (see Corollary~1 to Theorem~2.1, \cite{ParRanVar-1967ab} 
    or Theorem~5.1, \cite{Kum-2010aa}).
\end{proof}

\blmma \label{lm:Connb} 
Let $\lieg$ be a simple complex Lie algebra and let 
$\Lambda,\Gamma\in P_{+}$.
Let $B:=B(\Lambda)\otimes B(-\omega_{0}.\Gamma)$ be the product crystal of 
the tensor product module $V:=V(\Lambda)\otimes V(-\omega_ {0}.\Gamma)$. 
Let $b:=b^{\Lambda}_{\Lambda}\otimes b^{-\omega_{0}.\Gamma}_ {-\Gamma}$ 
be the vector in $B$ of weight $\Lambda-\Gamma$. Then,
 the connected component in $B$ containing $b$, say $\conn(b)$, has the
 highest weight $\dom(\Lambda-\Gamma)$ or, in other words, $\conn(b)$ 
 corresponds to the PRV component $V(\dom(\Lambda-\Gamma))$ in $V$. 
\elmma
 To prove this, one has to take an alternate approach to Kashiwara 
 crystals, namely, Littelmann path 
 crystals (Corollary 6.4.27, \cite{Jos-1995aa}).
 Then the proof of the above Lemma follows as a byproduct of Littelmann's 
 proof of the PRV conjecture using path crystals. For details, we refer 
 to sections~6.4.17 and 6.4.18 in \cite{Jos-1995aa} and section~7 
 in $\cite{Lit-1994aa}$.

\bthm
Let $\lieg$ be a simple complex Lie algebra.
For any dominant weight $\Omega\neq 0$, choose a lower global basis 
vector $v^{\Omega}_{j}$ of $V(\Omega)$ with $\wt(v^{\Omega}_{j})\neq 0$. 
Then the following are equivalent.
\begin{enumerate}
\item  
There are dominant weights $\Lambda,\Gamma$, and a $\utg$-module morphism
\begin{center}
 $T: V(\Lambda)\otimes V(-\omega_{0}.\Gamma)\to V(\Omega)$ 
\end{center}
        such that $T( L(\Lambda)\otimes L(-\omega_{0}.\Gamma))=L(\Omega)$ and 
        $T(v^{\Lambda}_{\Lambda}\otimes v^{-\omega_{0}.\Gamma}_{-\Gamma})=v^{\Omega}_{j}$. 
\item 
$\wt(v^{\Omega}_{j})=\omega. \Omega$ for some $\omega\in W$.
\end{enumerate}
\ethm
\begin{proof}
Denote the tensor product module $V(\Lambda)\otimes V(-\omega_{0}.\Gamma)$ by 
$V$ and the product crystal basis by $(L(V),B(V))$. 
Let us first prove the implication $(2)\implies (1)$. 
Choose $\Lambda$ and $\Gamma$ as in Lemma~\ref{lm:PRV}. 
Then, using the uniqueness of the crystal basis, we can find an 
$\utg$-module isomorphism
\[
S: V\to V(\Omega)\bigoplus(\bigoplus_{\Omega^{1}\neq \Omega} V(\Omega^{1}))
\] 
that sends $(L(V),B(V))$ to $(L(\Omega)\bigoplus(\bigoplus_{\Omega^
 {1}\neq \Omega}L(\Omega^{1}), B(\Omega)\bigsqcup(\bigsqcup_{\Omega^
 {1}\neq \Omega} B(\Omega^{1})))$. Denote $v^{\Lambda}_{\Lambda}\otimes v^
 {-\omega_{0}.\Gamma}_{-\Gamma}$ by $v$ and its projection to $B(V)$ by $b$. Write
 $S(v)=v_{\Omega}\bigoplus(\bigoplus_{\Omega^{1}\neq \Omega}\; v_{\Omega^
 {1}})$.

We now claim that $v_{\Omega}\notin t.L(\Omega)$. 
To prove this, assume, on the contrary, that
$v_{\Omega}\in t.L(\Omega)$. As $S$ is an isomorphism of crystal bases,
there is exactly one $\Omega^{1}(\neq \Omega)$ for which $v_{\Omega^
{1}}\notin t.L(\Omega^{1})$. So the induced map $\tilde{S}$, which is a
bijection from $B(V)$ to 
$B(\Omega)\bigsqcup(\bigsqcup_{\Omega^{1}\neq \Omega} B(\Omega^{1})$, 
sends $b$ to $v_{\Omega^{1}}+t.L(\Omega^{1})\in B(\Omega^{1})$. 
Now it follows that the connected component $\conn(b)$ of $B(V)$ 
is isomorphic to $B(\Omega^{1})$ as crystals via $\tilde{S}|_{\conn(b)}$. 
On the other hand, $\conn(b)$ corresponds to the PRV component 
$V(\dom(\Lambda-\Gamma))$ in $V$. But $\Lambda$ and $\Gamma$ are chosen in a
way so that $\dom(\Lambda-\Gamma)=\Omega$, and hence $\conn(b)$ has the
highest weight $\Omega$. But  this forces us to conclude that 
$\Omega=\Omega^{1}$, which is a contradiction.
Therefore, $\tilde{S}(b)=v_{\Omega}+t.L(\Omega)\in B(\Omega)$. 
The crystal $B(\Omega)$ has exactly one vector of weight 
$\Lambda-\Gamma=\omega.\Omega$, which is $v^{\Omega}_{\omega.\Omega}+t.L(\Omega)$. 
But $v_{\Omega}+t.L(\Omega)\in B(\Omega)$ is also of weight $\omega.\Omega$. 
Now it follows that $v_{\Omega}=c.v^{\Omega}_{\omega.\Omega}$ where $c$ 
is a unit in $A_{0}$. Thus the map $T:=(\frac{1}{c} I\oplus 0)\circ S$ 
will have all the desired properties.

Next we proceed to prove the converse implication $(1)\implies (2)$. 
In this case, $T$ induces a map
\[
\tilde{T}:L(V)/t.L(V)\to L(\Omega)/t.L(\Omega).
\]
The map $\tilde{T}$ commutes with Kashiwara operators 
$\tilde{E}_{i}, \tilde{F}_{i}$ , and $\tilde{T}(b)=b^{\Omega}_{j}$, 
therefore we have $\tilde{T}(\conn(b))\subseteq B(\Omega)\bigsqcup \{0\}$. 
Denote the restriction of $\tilde{T}$ to $\conn(b)$ by $\tilde{T}_{b}$. 
Since $T$ is a $\utg$-module morphism, it preserves weights. 
Hence  for any element $b_{1}\in \conn(b)$, one has
$\wt(\tilde{T}_{b}(b_1))=\wt(b_1)$ if $\tilde{T}_{b}(b_1)\neq 0$. 
Let $b_{h}$ be the highest weight vector of $\conn(b)$. Then it is of weight 
$\dom(\Lambda-\Gamma)$ and we must have $\tilde{T}_{b}(b_{h})\neq 0$, and
$\tilde{E}_{i}. \tilde{T}_{b}(b_{h})=0$ for all $i$. From this, it follows
that $0\neq\tilde{T}_{b}(b_{h})\in \mathbb{Q}.b^{\Omega}_{\Omega}$, and hence
$\wt(v^{\Omega}_{j})=\Lambda-\Gamma\in W.\Omega$.
\end{proof}
As an immediate consequence, we have the following corollary.
\bcrlre  \label{cr:Triangle}
 Let $\Omega\in P_{+}$, $\omega\in W$ and $\wt(v^{\Omega}_{j})=\omega.\Omega$.
 Then for any $1\leq i\leq\dim\Omega$, one has
$C^{\Omega}_{v^{\Omega}_{i},v^{\Omega}_{j}}\in R^{A_{0}}_{+}.R^{A_{0}}_{-}$. 
\ecrlre
We will now describe the crystal lattice $\otgazero$ 
as an $A_{0}$-algebra generated by a certain collection of matrix elements.
\bppsn
\label{pr:otgazero-as-algebra}
Let $\lieg$ be a simple complex Lie algebra and $G$ be a connected simply connected  Lie group with  $\mathrm{Lie}\,G=\lieg$.
Let $F$ be a set of highest weights such that for any highest weight
$\Lambda$, the highest weight module $V(\Lambda)$ is a submodule 
of a tensor product
$\otimes_{k=1}^{r} V(\Omega_{k})$ where $\Omega_{k}\in F$ for $1\leq k\leq r$.
Then $\otgazero$ is the $A_{0}$-algebra generated by the matrix elements
\begin{equation}\label{eq:gen-set}
\Bigl\{C^{\Omega}_{i,j}: \Omega\in F, \, 1\leq i,j\leq\dim\Omega\Bigr\}.
\end{equation}
\eppsn
\bprf
Since $\otgazero$ is an $A_{0}$-algebra, in view of (\ref{eq:otgazero}), 
it is enough to show that any matrix element
$C^{\Lambda}_{v_{i}^{\Lambda},v_{j}^{\Lambda}}$
belongs to the $A_{0}$-algebra generated by (\ref{eq:gen-set}).
Assume $V(\Lambda)$ is a direct summand in the tensor product
$\otimes_{k=1}^{r} V(\Omega_{k})$, where $\Omega_{k}\in F$ for $1\leq k\leq r$.
Then there is an $U_{t}(\mathfrak{g})$-module morphism $T$ from  
$\otimes_{k=1}^{r} V(\Omega_{k})$ onto $V(\Lambda)$ such that 
it carries the crystal lattice of the product module onto the crystal lattice
$L(\Lambda)$ in $V(\Lambda)$. Therefore, there exists 
an element 
$w:=\sum_{i}f_{i}(t) v^{\Omega_{1}}_{k_{i,1}} \otimes 
      \ldots \otimes v^{\Omega_{r}}_{k_{i,r}}$, 
where $f_{i}(t)\in A_{0}$ for all $i$,
such that 
$Tw=v^{\Lambda}_{j}$. Then for any $a\in \utg$,
\begin{IEEEeqnarray*}{rCl}
C^{\Lambda}_{v_{i}^{\Lambda},v_{j}^{\Lambda}}(a)
  &=& \langle(v^{\Lambda}_{i})^{*}, a Tw\rangle\\
  &=& \langle (v^{\Lambda}_{i})^{*}, Ta w\rangle\\
  &=& \langle T^{tr}(v^{\Lambda}_{i})^{*}, a w\rangle.
\end{IEEEeqnarray*}
Since $T^{tr}(v^{\Lambda}_{i})^{*}$ is a $A_{0}$-linear combination
of functionals of the form 
$(v^{\Omega_{1}}_{k_{1}})^{*} \otimes 
      \ldots \otimes (v^{\Omega_{r}}_{k_{r}})^{*}$,
the result follows.
\eprf

\bthm \label{th:Minuscule} 
Let $\mathfrak{g}$ be a complex simple Lie algebra of type
$A_{n}$, $B_{n}$, $C_{n}$, $D_{n}$, $E_6$ or $E_{7}$, and
$G$ be the connected simply connected Lie group such that
$\mathfrak{g}=\mathrm{Lie}\,G$. Then the crystal lattice 
$\mathcal{O}^{A_{0}}_{t}(G)$ is the $A_{0}$-algebra generated 
by $R^{A_{0}}_{+}$ and $R^{A_{0}}_{-}$.
\ethm
\bprf
If $\mathfrak{g}$ is of  type $A_{n}$, $B_{n}$, $C_{n}$, $E_6$ or $E_{7}$,
then one can take $F$ in Proposition~\ref{pr:otgazero-as-algebra} to be
a singleton set consisting of one minuscule dominant weight and if $\mathfrak{g}$ 
is of  type $D_{n}$, one can take $F$ to be a doubleton consisting of 
two minuscule dominant weights (for details, see section~I.2.1, \cite{Ros-1990aa}). By Corollary~\ref{cr:Triangle}, it follows that 
$C^{\Lambda}_{v^{\Lambda}_{i},v^{\Lambda}_{j}}\in R^{A_{0}}_{+}R^{A_{0}}_{-}$
for all $1\leq i,j\leq\dim\Lambda$.
Therefore by Proposition~\ref{pr:otgazero-as-algebra}, it follows that
$\mathcal{O}^{A_{0}}_{t}(G)$ is the $A_{0}$-algebra generated 
by $R^{A_{0}}_{+}$ and $R^{A_{0}}_{-}$.
\eprf

If $\mathfrak{g}$ is of  type $E_{8}$, $F_{4}$ or $G_{2}$, one can still take
$F$ to be singleton, but there are no minuscule dominant weights in these
cases (for a complete list of minuscule and quasi-minuscule representations, 
we refer to Appendix~A.2.3, \cite{LakRag-2008aa}). In each of these cases, 
one has a quasi-minuscule dominant weight as a tensor generator 
(section~I.2.1, \cite{Ros-1990aa}), which has a nonzero weight space 
corresponding to the weight 0. The weight vectors of weight 0 do not allow 
the above argument to go through, as the following proposition shows.
\bppsn
Let $\lieg$ be a simple complex Lie algebra and let $\Lambda,\Omega\in P_{+}$,
with $\Omega\neq 0$.
Let $v^{\Omega}_{j}$ be a nonzero vector in $V(\Omega)$ of weight 0. Then there does not exist any $\utg$-module morphism
\[
T: V(\Lambda)\otimes V(\Lambda)^{*} \to V(\Omega)
\]
such that $T(L(\Lambda)\otimes L(-\omega_{0}.\Lambda))=L(\Omega)$ and 
$T(v_{\Lambda}^{\Lambda}\otimes v_{-\Lambda}^{-\omega_{0}.\Lambda})=v^{\Omega}_{j}$.
\eppsn
\begin{proof}
Assume that such a $T$ exists. Then $T$ induces a non-zero 
$\mathbb{Q}$-linear map
\[
T_{0}: L(V)/t.L(V) \to L(\Omega)/t.L(\Omega)
\]
that sends $b:=b^{\Lambda}_{\Lambda}\otimes b^{-\omega_{0}.\Lambda}_{-\Lambda}$ 
to $b^{\Omega}_{j}$. Since the crystal of $B(\Omega)$ is connected, there exist
$i_{1},i_{2},\ldots,i_{k}$ such that 
$b^{\Omega}_{\Omega}= \tilde{E}_{i_1}\tilde{E}_{i_2}\ldots\tilde{E}_{i_k}b^{\Omega}_{j}$. 
However, from the tensor product rule of the product crystal, it follows
that the connected component of the product crystal containing $b$ is $\{b\}$.
Therefore, 
\[
T_{0}(0)=T_{0}(\tilde{E}_{i_1}\tilde{E}_{i_2}\ldots\tilde{E}_{i_k}b)
=b^{\Omega}_{\Omega},
\]
which gives us a contradiction.
\end{proof}

\subsection{The inclusion $\otgazero\subseteq \otkazero$}
Let us first look at how the matrix elements behave under the $*$-operation.
In this subsection, we will denote $\dim\Lambda$ by $m$.
\blmma
Let $\Lambda=\sum_{i}m_{i}\varpi_{i}\geq 
      -\omega_{0}\Lambda =\sum_{i}n_{i}\varpi_{i}$
in the lexicographic ordering. Let $\mu$ and $\nu$ be the 
weights of $v^{\Lambda}_{r}$ and $v^{\Lambda}_{s}$
respectively and let $w^{\Gamma}_{j}$'s be as in 
Proposition~\ref{ppsn:dual-lattice}. Then
\begin{IEEEeqnarray}{rCl}
\label{eq:star*}
(C_{r, s}^{\Lambda})^{*}
  &=& 
  c_{\Lambda} t^{(\mu - \nu, \rho)} 
  C^{-\omega_{0}\Lambda}_{w^{-\omega_{0}\Lambda}_{m-r+1}, 
     w^{-\omega_{0}\Lambda}_{m-s+1}}.
\end{IEEEeqnarray}
\elmma
\bprf
Let $a\in\utg$.
Then one has
\begin{IEEEeqnarray*}{rCl}
 \left(C^{\Lambda}_{r,s}\right)^{*}(a)
   &=& 
   C^{\Lambda}_{r,s}((Sa)^{*})\\
   &=&
   (v^{\Lambda}_{r},(Sa)^{*}v^{\Lambda}_{s})\\
   &=&
   ((Sa)v^{\Lambda}_{r},v^{\Lambda}_{s})\\
   &=&
   \langle (Sa)v^{\Lambda}_{r}, (v^{\Lambda}_{s})^{*}\rangle\\
   &=&
   \langle v^{\Lambda}_{r}, a\,(v^{\Lambda}_{s})^{*}\rangle.
\end{IEEEeqnarray*}
For any $f\in V(\Lambda)^{*}$ and a $v\in V(\Lambda)$ of weight $\mu$, 
we have
\begin{IEEEeqnarray*}{rCl}
 \langle v, f\rangle
   &=& 
   t^{(\mu,2\rho)}\langle K^{-2\rho}v, f\rangle\\
   &=&
   t^{(\mu,2\rho)}\left(K^{-2\rho}v, f_{*}\right)\\
   &=&
   t^{(\mu,2\rho)}(v^{\Lambda}_{m}, v^{\Lambda}_{m}) 
     t^{(-\Lambda, 2\rho)}\left(v^{*}, (f_{*})^{*}\right)\\
   &=&
   c_{\Lambda} t^{(\mu-\Lambda,2\rho)}
       \left(v^{*}, f \right).
\end{IEEEeqnarray*}
Thus we obtain
\begin{IEEEeqnarray*}{rCl}
 \left(C^{\Lambda}_{r,s}\right)^{*}(a)
   &=& 
   \langle v^{\Lambda}_{r}, a\,(v^{\Lambda}_{s})^{*}\rangle\\
   &=&
   c_{\Lambda}t^{(\mu-\Lambda,2\rho)}
    \left( (v^{\Lambda}_{r})^{*}, a\,(v^{\Lambda}_{s})^{*}\right)\\
   &=&
   c_{\Lambda}t^{(\mu-\nu,\rho)}
    \left( w^{-\omega_{0}\Lambda}_{m-r+1},  
       a\,(w^{-\omega_{0}\Lambda}_{m-s+1})\right)\\
   &=&
   c_{\Lambda}t^{(\mu-\nu,\rho)}
   C^{-\omega_{0}\Lambda}_{w^{-\omega_{0}\Lambda}_{m-r+1}, 
           w^{-\omega_{0}\Lambda}_{m-s+1}}(a),
\end{IEEEeqnarray*}
which gives us (\ref{eq:star*}).
\eprf

\brmrk
\label{rm:star-dual}
Assume  all the weight spaces of $V(\Lambda)$ are one dimensional.
Then $w^{-\omega_{0}\Lambda}_{k}=c. v^{-\omega_{0}\Lambda}_{k}$
for some $c\in \pm 1+tA_{0}$. Hence there exists an $e\in \pm 1+tA_{0}$
such that
\begin{IEEEeqnarray}{rCl}
\label{eq:star-onedim}
(C^{\Lambda}_{\mu, \nu})^{*}
  &=& 
  e t^{(\mu - \nu, \rho)} 
  C^{-\omega_{0}\Lambda}_{v^{-\omega_{0}\Lambda}_{-\mu}, 
     v^{-\omega_{0}\Lambda}_{-\nu}}.
\end{IEEEeqnarray}
Note that the quantity $(\mu - \nu, \rho)$ is an integer.
\ermrk

We have seen in Subsection~3.1 that the $A_{0}$-algebra 
$\mathcal{O}^{A_{0}}_t(G/N^{+})$ is the $A_{0}$-span
$R^{A_{0}}_{+}$ of the matrix elements 
$C^{\Lambda}_{i,1}\equiv C^{\Lambda}_{v^{\Lambda}_{i}, v^{\Lambda}_{\Lambda}}$. 
We will next give a similar description of the algebra 
$\mathcal{O}^{A_{0}}_t(G/N^{-})$ in terms of the matrix
elements 
$C^{\Lambda}_{i,m}\equiv C^{\Lambda}_{v^{\Lambda}_{i}, 
     v^{\Lambda}_{\omega_{0}\Lambda}}$.
\bppsn
\label{ppsn:Rplus-star}
Define
\begin{IEEEeqnarray*}{rCl}
\tilde{R}_-^{A_0} &:=& 
 \mathrm{Span}_{A_0} \{t^{(\omega_{0}\Lambda - \wt(v^{\Lambda}_{i}), \rho)} 
  C^{\Lambda}_{i, m} : \Lambda \in P_+, 1\leq i\leq m\}.
\end{IEEEeqnarray*}
Then one has
\begin{IEEEeqnarray}{rClrCl}
(R_+^{A_0})^{*}&=& \tilde{R}_-^{A_0},\qquad &
(R_{-}^{A_0})^{*}&\subseteq & R_{+}^{A_0}.
\label{eq:rstars}
\end{IEEEeqnarray}
\eppsn
\bprf
It follows from (\ref{eq:star*}) that if $v^{\Lambda}_{i}$ is
of weight $\mu$, then for
$\Lambda\geq -\omega_{0}\Lambda$ (in lexicographic ordering),
one has
\begin{IEEEeqnarray}{rCl}
\label{eq:star1}
(C^{\Lambda}_{v^{\Lambda}_{i}, v^{\Lambda}_{1}})^{*}
  &=& 
  c_{\Lambda} t^{(\mu - \Lambda, \rho)} 
C^{-\omega_{0}\Lambda}_{w^{-\omega_{0}\Lambda}_{m-i+1}, 
           w^{-\omega_{0}\Lambda}_{m}},\\
\label{eq:star3}
(C_{v^{\Lambda}_{i}, v^{\Lambda}_{\omega_{0}\Lambda}}^{\Lambda})^{*}
  &=& 
  c_{\Lambda} t^{(\mu - \omega_{0}\Lambda, \rho)} 
C^{-\omega_{0}\Lambda}_{w^{-\omega_{0}\Lambda}_{m-i+1}, 
     w^{-\omega_{0}\Lambda}_{1}},\\
\label{eq:star2}
\left(C^{-\omega_{0}\Lambda}_{w^{-\omega_{0}\Lambda}_{m-i+1}, 
     w^{-\omega_{0}\Lambda}_{m}}\right)^{*}
    &=&
    c_{\Lambda}^{-1} t^{(\Lambda - \mu, \rho)}
    C^{\Lambda}_{v^{\Lambda}_{i}, v^{\Lambda}_{1}},\\
\label{eq:star4}
\left(C^{-\omega_{0}\Lambda}_{w^{-\omega_{0}\Lambda}_{m-i+1}, 
   w^{-\omega_{0}\Lambda}_{1}}\right)^{*}
    &=&
    c_{\Lambda}^{-1} t^{(\omega_{0}\Lambda - \mu, \rho)} 
    C^{\Lambda}_{v^{\Lambda}_{i}, v^{\Lambda}_{\omega_{0}\Lambda}}.
\end{IEEEeqnarray}
Let us now consider highest weights that are self-dual and those that are not
separately.
\paragraph{Case I: $\Lambda\neq -\omega_{0}\Lambda$.}\mbox{}\newline
Assume $\Lambda > -\omega_{0}\Lambda$. Then from 
Proposition~\ref{ppsn:dual-lattice}, Remark~\ref{rm:glb-dual}
and (\ref{eq:star1}), it follows that
$(C^{\Lambda}_{r,1})^{*}\in \tilde{R}^{A_{0}}_{-}$.
Similarly, using Proposition~\ref{ppsn:dual-lattice},
Remark~\ref{rm:glb-dual}  and (\ref{eq:star4}),  we obtain
$(C^{-\omega_{0}\Lambda}_{r,1})^{*}\in \tilde{R}^{A_{0}}_{-}$.
Therefore in this case, we have $(R_+^{A_0})^{*}\subseteq \tilde{R}_-^{A_0}$.
\paragraph{Case II: $\Lambda = -\omega_{0}\Lambda$.}\mbox{}\newline 
Let $T$ be the $\utg$-module isomorphism from $V(\Lambda)$ to 
$V(\Lambda)^{*}$ as in (\ref{eq:iso-to-dual}).
Then $T(v^{\Lambda}_1) = w^{-\omega_{0}\Lambda}_{1}$,
$(Tx, Ty) = (x, y)$ for all $x,y \in V(\Lambda)$,
and the restriction of $T$ to the crystal lattice of 
$V(\Lambda)$ is an $A_0$-linear isomorphism from $L(\Lambda)$ 
onto $L(-\omega_{0}\Lambda)$.  Therefore, 
(\ref{eq:star1}) can be further written as
\[
(C_{r,1}^{\Lambda})^{*}= 
c_{\Lambda} t^{(\mu - \Lambda, \rho)} 
C^{-\omega_{0}\Lambda}_{w^{-\omega_{0}\Lambda}_{m-r+1}, w^{-\omega_{0}\Lambda}_{m}}
= c_{\Lambda} t^{(\mu - \Lambda, \rho)} 
C^{\Lambda}_{T^{-1}(w^{-\omega_{0}\Lambda}_{m-r+1}), 
     T^{-1}(w^{-\omega_{0}\Lambda}_{m})},
\]
where $\mu = \wt(v^{\Lambda}_r)$.
Note that $T^{-1}(w^{-\omega_{0}\Lambda}_{m-r+1})$ and 
$T^{-1}(w^{-\omega_{0}\Lambda}_{m})$ have weights $-\mu$ and $-\Lambda$ 
respectively, and they lie in the weight spaces $L(\Lambda)_{-\mu}$ 
and $L(\Lambda)_{-\Lambda}$ of $L(\Lambda)$. Hence  
$T^{-1}(w^{-\omega_{0}\Lambda}_{m-r+1})$ can be expressed 
as $A_0$-linear combination of lower global basis vectors of 
$V(\Lambda)$ having weight $-\mu$. Similarly, 
$T^{-1}(w^{-\omega_{0}\Lambda}_{m})$ is $A_0$-multiple of $v^{\Lambda}_{m}$. 
Now, using the identity $\Lambda = -\omega_{0}\Lambda$, we can 
write $t^{(\mu - \Lambda, \rho)} = t^{(\omega_{0}\Lambda - (-\mu), \rho)}$, 
and thus we have
\[
(C_{r,1}^{\Lambda})^{*}= \sum_k d_k \cdot t^{(\omega_{0}\Lambda -(-\mu), \rho)} 
  C_{v^{\Lambda}_{j_k}, v^{\Lambda}_m}^{\Lambda},
\] 
where $d_k \in A_0$, and $\mathrm{span}_{A_0} \{v^{\Lambda}_{j_k} \in G
(\Lambda): \; k\} = L(\Lambda)_{-\mu}$. 

Combining both the cases, we finally
have $(R^{A_{0}}_{+})^{*}\subseteq \tilde{R}^{A_{0}}_{-}$. The opposite
inclusion follows similarly using (\ref{eq:star4}) and 
Proposition~\ref{ppsn:dual-lattice}. 
Thus we have the first equality in(\ref{eq:rstars}).
Since $R^{A_{0}}_{-}\subseteq \tilde{R}_-^{A_0}$, the inclusion
$(R_{-}^{A_0})^{*}\subseteq  R_{+}^{A_0}$ follows immediately.
\eprf
 
\bcrlre
\label{cr:rtilde}
The space $\tilde{R}_-^{A_0}$ is an $A_{0}$-subalgebra of
$\otg$ and one has $\mathcal{O}^{A_{0}}_t(G/N^{-})=\tilde{R}_-^{A_0}$.
\ecrlre
\bprf
Since $R_{+}^{A_0}$ is an $A_{0}$-algebra, it follows that so is
$\tilde{R}_-^{A_0}$. It now follows from Proposition~3.3 \cite{MatYun-2023aa},
Proposition~\ref{ppsn:rplus} and Proposition~\ref{ppsn:Rplus-star} above
that $\mathcal{O}^{A_{0}}_t(G/N^{-})=\tilde{R}_-^{A_0}$.
\eprf
Recall (Definition~3.2, \cite{MatYun-2023aa}) that $\otkazero$ is the 
$A_{0}$-subalgebra of $\otg$ generated by $\mathcal{O}_{t}^{A_{0}}(G/N^{+})$
and $\mathcal{O}_{t}^{A_{0}}(G/N^{-})$. Using Proposition~\ref{ppsn:Rplus-star} and
Theorem~\ref{th:Minuscule}, the following inclusion which was conjectured by 
Matassa \& Yuncken (see Remark~3.4, \cite{MatYun-2023aa}) is now immediate when
$\mathfrak{g}$ is of type  $A_{n}$, $B_{n}$, $C_{n}$, $D_{n}$, $E_6$ or $E_{7}$.
\bcrlre
\label{cr:M-Yconj}
    $\otgazero\subseteq \mathcal{O}^{A_{0}}_{t}(K)$.
\ecrlre
It follows from the above that the $\ast$-algebra $\tlotgazero$ is contained in
$\otkazero$.  Thus we have $\otkazero=\tlotgazero$, i.e.\ $\otkazero$ is the 
$*$-algebra over $A_{0}$ generated by $\otgazero$. This is a key observation 
that we will use in the next section as well as in Section~5 while comparing two 
different approaches to crystal limits for the type $A_{n}$ case.

\section{Matassa-Yuncken crystallization} 

\subsection{The crystallized algebra} 
Let us start by recalling the Matassa-Yuncken definition of the
crystallized algebra $C(K_{0})$, which we will denote by $C_{MY}$.
They first define the space  $\otkazero$ which we have encountered
in Sections~4 and 5.
Then for each $a\in \otkazero$, they produce an operator $\pi_{0}(a)$
in $\mathcal{L}(\mathcal{H}_{Soi})$ by using the following mechanism.

For each $i\in\{1,2,\ldots,n\}$, let $\phi^{}_{t,i}$ denote the map
from $\otg$ onto   $\otisltc$ dual to the inclusion map 
$\phi^{U}_{t,i}:\utisltc\to \utg$ given by
\[
E\mapsto E_{i},\quad F\mapsto F_{i},\quad K\mapsto K_{i}.
\]
Simple computations tell us that the map $\phi^{}_{t,i}$ is 
a $\ast$-homomorphism  from $\otg$ onto $\otisltc$.
This map $\phi^{}_{t,i}$ carries the space 
$\otgazero$ to $\otisltcazero$ and consequently, the space
$\otkazero$ to $\otisutwoazero$. 
They now use what they call
a `patially defined evaluation map' $ev_{q}$ from $\otsutwoazero$
to $\oqsutwo$ to come down to the space $\oqisutwo$. We will denote this
map by $\vartheta_{q}$ below.

Let $\pi_{q}:\mathcal{O}_{q}(SU(2))\to \mathcal{L}(\ell^{2}(\mathbb{N}))$
and $\chi_{q}:\mathcal{O}_{q}(SU(2))\to \mathcal{L}(\ell^{2}(\mathbb{Z}))$
be the representations given by
\[
\pi_{q}(C^{\one}_{i,j})=
    \begin{cases}
      S\sqrt{I-q^{2N}} &\text{if }(i,j)=(1,1),\cr
      \sqrt{I-q^{2N}}S^{*} & \text{if } (i,j)=(2,2),\cr
      -q^{N+1} &  \text{if } (i,j)=(1,2),\cr
      q^{N}   &  \text{if }(i,j)=(2,1),
    \end{cases}
\qquad
\chi_{q}(C^{\one}_{i,j})=
    \begin{cases}
      S^{*} &\text{if }(i,j)=(1,1),\cr
      S & \text{if } (i,j)=(2,2),\cr
      0 &  \text{otherwise},
    \end{cases}
\]
where we denote by $S$ the left shift $e_{j}\mapsto e_{j-1}$
both on $\ell^{2}(\mathbb{N})$ as well as on $\ell^{2}(\mathbb{Z})$ 
and by $N$ the operator $e_{j}\mapsto je_{j}$ so that
$q^{f(N)}$ is the operator $e_{j}\mapsto q^{f(j)}e_{j}$.
Now
take the longest word
$\omega_{0}=s_{i_{1}}s_{i_{2}}\ldots s_{i_{k}}$ 
of the Weyl group. Then for $a\in\otkazero$,
\begin{IEEEeqnarray}{rCl}
\tilde{\pi}_{q}(a) 
 &=&\lim_{q\to 0+}\Bigl(
   (\pi_{q}\circ\vartheta_{q}\circ\phi^{}_{t,i_{1}})\otimes\cdots\otimes
   (\pi_{q}\circ\vartheta_{q}\circ\phi^{}_{t,i_{k}})\otimes{}\Bigr.\nonumber\\
  && \hspace{1em}\Bigl.(\chi_{q}\circ\vartheta_{q}\circ\phi^{}_{t,1})\otimes
   (\chi_{q}\circ\vartheta_{q}\circ\phi^{}_{t,2})\otimes\cdots\otimes
   (\chi_{q}\circ\vartheta_{q}\circ\phi^{}_{t,n})\Bigr)\Delta^{(k+n-1)}(a),
     \label{eq:tildepiq}\\
  \pi_{0}(a)  &=& \lim_{q\to 0+} \tilde{\pi}_{q}(a).
\end{IEEEeqnarray}
The crystallized algebra $C_{MY}$ is then defined to be the $C^{*}$-subalgebra
of $\mathcal{L}(\mathcal{H}_{Soi})$ generated by the $\pi_{0}(a)$'s for
$a\in \otkazero$. 

\brmrk
Matassa \& Yuncken  describes the torus part 
$\chi:\oqsu\to \mathcal{L}(L^{2}(\mathbb{T}^{n}))
  \cong\mathcal{L}(\ell^{2}(\mathbb{Z}^{n}))$ 
of the Soibelman representation 
(see the discussion preceeding Definition 4.4, \cite{MatYun-2023aa}),
and define $\pi_{0}$ in terms of its `partially defined analogue on $\otk$'.
The $\chi$ described by them can easily be seen to be unitarily equivalent
to 
$\Bigl(\chi_{q}\otimes \chi_{q}\otimes\cdots\otimes
   \chi_{q}\Bigr)\Delta^{(n-1)}(\cdot)$.
The `partially defined analogue', however, is not clearly specified
in \cite{MatYun-2023aa}.
We take this `partially defined analogue' to be
$\Bigl.(\chi_{q}\circ\vartheta_{q}\circ\phi^{}_{t,1})\otimes
   (\chi_{q}\circ\vartheta_{q}\circ\phi^{}_{t,2})\otimes\cdots\otimes
   (\chi_{q}\circ\vartheta_{q}\circ\phi^{}_{t,n})\Bigr)\Delta^{(n-1)}(\cdot)$,
in line with the partially defined analogue in the non-torus part
that they work with, which is projecting down to $\otsutwo$ first and 
then using the evaluation map $\vartheta_{q}$.
\ermrk

We will next focus on the evaluation map $\vartheta_{q}$.
This map is used in Theorem~4.2 and Corollary~4.3, \cite{MatYun-2023aa} to define
$\tilde{\pi}_{q}(a)$ for $a\in\otsltcazero$ and $a\in \otsutwoazero$ respectively,
and these two results are central to the definition of their crystallized 
algebra. However, the map $\vartheta_{q}$ is not defined clearly in their paper.

In the next subsection we will define, for any given $a\in\otg$,
an element $\vartheta_{q}(a)\in\oqg$ for each $q$ in a small interval
of the form $(0,\delta_{a})$ and prove that 
$\vartheta_{q}:\tlotgazero\to \oqk$ has the required
multiplicative and $\ast$-preserving properties. We will also show that 
instead of restricting to $\otsutwoazero$ and then using the specialization map and 
further using the representations of $\oqsutwo$ to construct $\tilde{\pi}_{q}(a)$
as has been done in \cite{MatYun-2023aa},
one could first specialize to $\oqk$ and then directly use the 
Soibelman representation of $\oqk$ to arrive at $\tilde{\pi}_{q}(a)$. Moreover,
this is valid not only for $a\in\otkazero$ but for all $a\in\tlotgazero$.

\subsection{The specialization maps} 
The $A$-form $\utga$ (also referred to as the 
restricted $A$-form or the Lusztig $A$-form) 
of $\utg$ is the $A$-subalgebra of $\utg$ generated by the elements 
\[
E_{i}^{(m)}:=E_{i}^{m}/[m]_{t^{d_{i}}}!, \quad
  F_{i}^{(m)}:=F_{i}^{m}/[m]_{t^{d_{i}}}!,\quad
    K_{i}, K_{i}^{-1},\qquad 1\leq i\leq n, \quad m\in\mathbb{N}.
\]
This is a Hopf-$\ast$ algebra over $A$.
Let $V(\Lambda)$ be a finite dimensional irreducible highest weight module
of $\utg$ with highest weight vector $v^{\Lambda}_{\Lambda}$. 
Then $V^{A}(\Lambda):=\utga v^{\Lambda}_{\Lambda}$ is an indecomposable
$\utga$-module.
For each $q\in (0,\infty)$, $q\neq 1$, one can then put a natural Hopf-$\ast$ structure on
$\mathbb{C}\otimes_{A}\utga$ over $\mathbb{C}\otimes_{A}A\cong\mathbb{C}$, 
where the tensor product over $A$ in both cases 
is via the action $t\mapsto q$ of $A$ on $\mathbb{C}$. 
One has an isomorphism of Hopf-$\ast$ algebras
\[
\uqg\cong \mathbb{C}\otimes_{A}\utga,
\]
and a resulting isomorphism of $\uqg$-modules
\[
V(\Lambda)\cong \mathbb{C}\otimes_{A}V^{A}(\Lambda).
\]
The map $\theta^{U}_{q}: \utga\hookrightarrow 
    \mathbb{C}\otimes_{A}\utga \cong \uqg$
is called the specialization map. We will denote the 
evaluation map
$p[t,t^{-1}]\mapsto 1\otimes_{A}p[t,t^{-1}]=p[q,q^{-1}]$
from $A$ to $\mathbb{C}$ by $e_{q}$.

At the dual level, the $A$-lattice $\otga$ of $\otg$ is defined to be
the subspace
$\{f\in\otg: \langle f,a\rangle\in A \text{ for all }a\in\utga\}$.
This is an $A$-subalgebra of $\otg$ and one can show that
the matrix elements $C^{\Lambda}_{(v^{\Lambda}_{i})^{*},v^{\Lambda}_{j} }$
as well as $C^{\Lambda}_{(v^{\Lambda}_{i})^{o},v^{\Lambda}_{j} }$
are elements of $\otga$.
There is a pairing between $\mathbb{C}\otimes_{A}\utga$ 
and $\mathbb{C}\otimes_{A}\otga$ given by
\[
\langle z\otimes_{A}f, w\otimes_{A} a\rangle
   := zw \,e_{q}(\langle f,  a\rangle) \in\mathbb{C}.
\]
This gives a Hopf-$\ast$ algebra isomorphism
\[
\oqg\cong \mathbb{C}\otimes_{A}\otga,
\]
and a resulting specialization map
$\theta_{q}: \otga\hookrightarrow 
    \mathbb{C}\otimes_{A}\otga \cong \oqg$.
For more details on $A$-forms and the specialization map, 
we refer the reader to \cite{Lus-1988aa}, \cite{Lus-1990ac}, \cite{Lus-2009aa}, 
Chapter~31, \cite{Lus-1993aa} and Chapter 9, \cite{ChaPre-1995aa}.

Note that the map $\theta_{q}$ is linear ($A$-linearlity being carried to
$\mathbb{C}$-linearity through the action $t\mapsto q$),
multiplicative, preserves the involution,  and commutes with the coproduct.
Take an $a\in\otg$. Since 
$\{C^{\Lambda}_{(v^{\Lambda}_{i})^{o},v^{\Lambda}_{j}}:
    \Lambda\in P_{+}, 1\leq i,j\leq\dim\Lambda\}$
is a $\mathbb{Q}(t)$-basis for $\otg$, the element $a$ is 
a finite sum of the form
$\sum_{\Lambda,i,j}f^{\Lambda}_{i,j}(t)
  C^{\Lambda}_{(v^{\Lambda}_{i})^{o},v^{\Lambda}_{j}}$,
where $f^{\Lambda}_{i,j}(t)$'s are from $\mathbb{Q}(t)$
and are unique. Hence there 
is a $\delta_{a}>0$ such that none of the 
$f^{\Lambda}_{i,j}(t)$'s have a pole in $(0,\delta_{a})$. 
Then for $q\in (0,\delta_{a})$,
define $\vartheta_{q}(a)$ to be 
$\sum_{\Lambda,i,j}f^{\Lambda}_{i,j}(q)
   \theta_{q}(C^{\Lambda}_{(v^{\Lambda}_{i})^{o},v^{\Lambda}_{j}})$.
Similarly, for $f(t)\in \mathbb{Q}(t)$ and for $q$ in 
a sufficiently small neighborhood $(0,\delta)$, 
define $e_{q}(f(t)):= f(q)$.

Note that the comultiplication map $\Delta$ of $\otg$
maps $\otgazero$ into $\otgazero\otimes\otgazero$.
Therefore it maps $\tlotgazero$ into $\tlotgazero\otimes\tlotgazero$.
For the rest of this subsection, let us denote the restriction
of the comultiplication map to $\tlotgazero$ by $\Delta_{t}$ and 
that of $\oqk$ by $\Delta_{q}$.
Let us follow a similar convention for the counit map also.
Using the fact that  $\theta_{q}$ is a Hopf-$\ast$ homomorphism, it 
follows that if $a$ and $b$ are two elements of
$\otg$, then one has, for a small enough $\delta>0$
and for $q\in (0,\delta)$,
\begin{IEEEeqnarray*}{rClrClrCl}
\vartheta_{q}(ab) &=& \vartheta_{q}(a)\vartheta_{q}(b),\qquad &
   \vartheta_{q}(a+b) &=& \vartheta_{q}(a)+\vartheta_{q}(b),\qquad
    & \vartheta_{q}(a^{*}) &=& (\vartheta_{q}(a))^{*},\\
(\vartheta_{q}\otimes \vartheta_{q})\Delta_{t}(a) 
    &=& \Delta_{q}(\vartheta_{q}(a)),
    & e_{q}\epsilon_{t}(a) &=& \epsilon_{q}\vartheta_{q}(a).
     &&&
\end{IEEEeqnarray*}

We will prove that the specialization map commutes with 
the restriction maps to copies of the subgroup $SL_{q}(2,\mathbb{C})$
and $SL_{t}(2,\mathbb{C})$.
\bppsn \label{pr:theta}
Let $\theta_{q}:\otka\to \oqk$ and $\vartheta^{}_{q}:\tlotgazero\to\oqk$ 
be the specialization maps introduced above.
For $1\leq i\leq n$, let $\phi^{}_{t,i}$ denote the homomorphism from 
$\otg$ onto $\otisltc$ defined earlier 
and let $\phi^{}_{q,i}$ denote the corresponding map from
$\oqg$ onto $\oqisltc$. Then $\phi^{}_{t,i}$ maps the lattices
$\otka$ and $\tlotgazero$ onto the respective lattices 
$\otisutwoa$ and $\tlotisltcazero=\otisutwoazero$ in $\otsltc$, $\phi^{}_{q,i}$ 
maps $\oqk$ onto $\oqisutwo$,
and the following diagrams commute:
\begin{center}
\begin{tikzcd}
\otka \arrow[r,
"\phi^{}_{t,i}"] \arrow[d, "\theta^{}_{q}"]
& \otisutwoa \arrow[d,
"\theta^{}_{q}" ] &&\tlotgazero \arrow[r,
"\phi^{}_{t,i}"] \arrow[d, "\vartheta^{}_{q}"]
& \otisutwoazero \arrow[d,
"\vartheta^{}_{q}" ] \\
\oqk \arrow[r, "\phi^{}_{q,i}"]
&  \mathcal{O}_{q_{i}}(SU(2))&&\oqk \arrow[r, "\phi^{}_{q,i}"]
&  \oqisutwo
\end{tikzcd}
\end{center}
\eppsn
\bprf
As noted before, we have 
$\mathbb{C}\otimes_{A} A\cong \mathbb{C}$
via the action $t\mapsto q$ and $e_{q}$
is the  evaluation map 
$f[t,t^{-1}]\mapsto 1\otimes_{A} f[t,t^{-1}]=f[q,q^{-1}]$
from $A$ into $\mathbb{C}\otimes_{A} A\cong \mathbb{C}$.
Then one has, for $a\in \utga$ and $f\in \otga$,
\[
\langle \theta^{U}_{q}(a), \theta_{q}(f)\rangle = 
  e_{q}(\langle a, f \rangle).
\]
Let $\phi_{t,i}^{U}$ denote the inclusion map
\[
E\mapsto E_{i},\quad F\mapsto F_{i},\quad K\mapsto K_{i}.
\]
from $\utisltca$ to $\utga$, and let $\phi_{q,i}^{U}$ denote the 
corresponding inclusion map from $\uqisltc$ to $\uqg$.
Observe that $\phi_{t,i}$ and $\phi_{q,i}$ are dual to the maps
$\phi_{t,i}^{U}$ and $\phi_{q,i}^{U}$ respectively.
It is straighforward to see that the following diagram commutes:
\begin{center}
\begin{tikzcd}
\utga  \arrow[d, "\theta^{U}_{q}"]
& \utisltca \arrow[l,"\phi^{U}_{t,i}"']\arrow[d,
"\theta^{U}_{q}" ] \\
\uqg 
&  \uqisltc\arrow[l, "\phi^{U}_{q,i}"']
\end{tikzcd}
\end{center}
Now take $a\in\utisltca$ and $f\in \otga$. Then we have
\begin{IEEEeqnarray*}{rCl}
\langle \theta^{U}_{q}(a),\theta_{q}\circ\phi_{t,i}(f) \rangle
  &=& 
  e_{q}(\langle a,\phi_{t,i}(f) \rangle)\\
  &=&
  e_{q}(\langle \phi^{U}_{t,i}(a), f \rangle)\\
  &=& 
\langle \theta^{U}_{q}\circ\phi_{t,i}^{U}(a),\theta_{q}(f) \rangle\\
  &=& 
\langle \phi_{q,i}^{U}\circ\theta^{U}_{q}(a),\theta_{q}(f) \rangle\\
  &=&
  \langle \theta^{U}_{q}(a),\phi_{q,i}\circ\theta_{q}(f) \rangle.
\end{IEEEeqnarray*}
Thus $\theta_{q}\circ\phi_{t,i}(f)=\phi_{q,i}\circ\theta_{q}(f)$ for all
$f\in \otga$, i.e.\ the first diagram commutes.

For commutativity of the second diagram, first take an $a\in \otgazero$.
Let $q>0$. Writing 
$a= \sum_{\Lambda,r,s}f^{\Lambda}_{r,s}(t)
  C^{\Lambda}_{(v^{\Lambda}_{r})^{o},v^{\Lambda}_{s}}$
and observing that 
$C^{\Lambda}_{(v^{\Lambda}_{r})^{o},v^{\Lambda}_{s}}\in\otka$,
it follows that both $\vartheta_{q}(a)$ and
$\vartheta_{q_{i}}(\phi_{t,i}(a))$ are defined 
if and only if $q\in(0,\delta_{a})$, where
\[
\delta_{a}=
\min\{z>0: z \text{ is a pole of } f^{\Lambda}_{r,s}(t)
     \text{ for some }\Lambda, r,s\}.
\]
The equality $\phi_{q,i}(\vartheta_{q}(a))=\vartheta_{q_{i}}(\phi_{t,i}(a))$
now follows from the commutativity of the first diagram. 
The same equality for $a\in\tlotgazero$ then follows immediately.
\eprf
Now we are in a position to have the following alternative description of the crystallized algebra of Matassa \& Yuncken.
\bthm
\label{thm:MY-alternate}
Let $\mathfrak{g}$ be a simple complex Lie algebra.
Then the Matassa-Yuncken crystallization of $C(K_{q})$ (respectively $\oqk$)
is the $C^{*}$-subalgebra (respectively $\ast$-subalgebra)
of $\mathcal{L}(\mathcal{H}_{Soi})$ generated by 
\[
\Bigl\{\lim_{q\to 0+}\psi^{(q)}_{Soi}\circ\vartheta_{q}^{}(a): a\in\otkazero\Bigr\}.
\]
\ethm
\bprf
Note that $\Delta(a)\in\tlotgazero\otimes\tlotgazero$ 
for all $a\in\tlotgazero$. 
Therefore $\tilde{\pi}_{q}(a)$ in (\ref{eq:tildepiq}) is defined
for $a\in\tlotgazero$ also, and for such an $a$, one has
{\allowdisplaybreaks
\begin{IEEEeqnarray*}{rCl}
 \tilde{\pi}_{q}(a)
   &=& 
\Bigl(
   (\pi_{q}\circ\vartheta_{q}^{}\circ\phi^{}_{t,i_{1}})\otimes\cdots\otimes
   (\pi_{q}\circ\vartheta_{q}^{}\circ\phi^{}_{t,i_{k}})\otimes{}\Bigr.\\
 &&\hspace{1em}\Bigl.(\chi_{q}\circ\vartheta_{q}^{}\circ\phi^{}_{t,1})\otimes
   (\chi_{q}\circ\vartheta_{q}^{}\circ\phi^{}_{t,2})\otimes\cdots\otimes
   (\chi_{q}\circ\vartheta_{q}^{}\circ\phi^{}_{t,n})\Bigr)
          \Delta^{(k+n-1)}(a)\\
   &=& 
\Bigl(
   (\pi_{q}\circ\phi^{}_{q,i_{1}}\circ\vartheta_{q}^{})\otimes\cdots\otimes
   (\pi_{q}\circ\phi^{}_{q,i_{k}}\circ\vartheta_{q}^{})\otimes{}\Bigr.\\
   &&\hspace{1em}\Bigl.(\chi_{q}\circ\phi^{}_{q,1}\circ\vartheta_{q}^{})\otimes
   (\chi_{q}\circ\phi^{}_{q,2}\circ\vartheta_{q}^{})\otimes\cdots\otimes
   (\chi_{q}\circ\phi^{}_{q,n}\circ\vartheta_{q}^{})\Bigr)
    \Delta^{(k+n-1)}(a)\\
   &=& 
\Bigl(
   (\pi_{q}\circ\phi^{}_{q,i_{1}})\otimes\cdots\otimes
   (\pi_{q}\circ\phi^{}_{q,i_{k}})\otimes{}\Bigr.\\
   &&\hspace{1em}\Bigl.(\chi_{q}\circ\phi^{}_{q,1})\otimes
   (\chi_{q}\circ\phi^{}_{q,2})\otimes\cdots\otimes
   (\chi_{q}\circ\phi^{}_{q,n})\Bigr)\circ
   \Bigl(\vartheta_{q}^{}\otimes\cdots\otimes\vartheta_{q}^{}\Bigr)
    \Delta^{(k+n-1)}(a)\\
   &=& 
\Bigl(
   (\pi_{q}\circ\phi^{}_{q,i_{1}})\otimes\cdots\otimes
   (\pi_{q}\circ\phi^{}_{q,i_{k}})\otimes{}\Bigr.\\
   &&\hspace{1em}\Bigl.(\chi_{q}\circ\phi^{}_{q,1})\otimes
   (\chi_{q}\circ\phi^{}_{q,2})\otimes\cdots\otimes
   (\chi_{q}\circ\phi^{}_{q,n})\Bigr)\circ\Delta^{(k+n-1)}
   (\vartheta_{q}^{}(a))\\
   &=&
   \psi^{(q)}_{Soi}\circ\vartheta_{q}^{}(a).
\end{IEEEeqnarray*}
}
Hence for all $a\in\otkazero\subseteq\tlotgazero$, we continue to have
\begin{IEEEeqnarray}{rCl}\label{eq:soi-evq}
\tilde{\pi}_{q}(a) &=& \psi^{(q)}_{Soi}\circ\vartheta_{q}^{}(a),
\end{IEEEeqnarray}
which gives us the required equality.
\eprf

\brmrk \label{rm:lim-tlotigazero}
\begin{enumerate}
	\item
 Since $\otisutwoazero=\tlotisltcazero$,
 it follows that 
$\lim_{q\to 0+}\psi^{(q)}_{Soi}\circ\vartheta_{q}^{}(a)$
exists for all $a\in\tlotgazero$.
    \item
    If $\mathfrak{g}$ is of one of the types listed in Theorem~\ref{th:Minuscule},
    then one has
\begin{IEEEeqnarray*}{rCl}
C_{MY} &= &
   \text{$C^{*}$-subalgebra of $\mathcal{L}(\mathcal{H}_{Soi})$ generated by }
 \Bigl\{\lim_{q\to 0+}\psi^{(q)}_{Soi}\circ\vartheta_{q}^{}(a): 
                a\in\otkazero\Bigr\}\\
  &=&
   \text{$C^{*}$-subalgebra of $\mathcal{L}(\mathcal{H}_{Soi})$ generated by }
  \Bigl\{\lim_{q\to 0+}\psi^{(q)}_{Soi}\circ\vartheta_{q}^{}(a): 
                a\in\tlotgazero\Bigr\}\\
  &=&
   \text{$C^{*}$-subalgebra of $\mathcal{L}(\mathcal{H}_{Soi})$ generated by }
 \Bigl\{\lim_{q\to 0+}\psi^{(q)}_{Soi}\circ\vartheta_{q}^{}(a): 
             a\in\otgazero\Bigr\}.
\end{IEEEeqnarray*}
\end{enumerate}
\ermrk

\section{Type $A_{n}$ case} 
We will prove in this section that the Matassa-Yuncken crystallized
algebra coincides with the cryatallized quantized function algebra
introduced by Giri \& Pal (\cite{GirPal-2022tv}) in the type $A_{n}$ case.
Let us start by describing 
the spaces $\otslncazero$ and $\otsuazero$ in terms of the
matrix elements $C^{\one}_{i,j}$. This will lead to a proof of 
equality of the two crystallizations mentioned above.
Observe that it follows from Proposition~\ref{pr:otgazero-as-algebra} 
that the crystal lattice $\otslncazero$ is the $A_{0}$-subalgebra of 
$\otslnc$ generated by $\Bigl\{C^{\one}_{i,j}: 1\leq i,j\leq n+1\Bigr\}$.

We will next express the space $\mathcal{O}_{t}^{A_{0}}(SU(n+1))$
as an algebra generated by the matrix elements $C^{\one}_{i,j}$ scaled 
appropriately.
\bthm
\label{th:otk}
The $A_{0}$-subalgebra $\mathcal{O}_{t}^{A_{0}}(SU(n+1))$ of 
$\otslnc$ is the $A_{0}$-algebra generated by 
$t^{\min\{i-j,0\}}C^{\one}_{i,j}$, $1\leq i,j\leq n$.
\ethm
\bprf
By Corollary~\ref{cr:M-Yconj}, $\mathcal{O}_{t}^{A_{0}}(SU(n+1))$
contains $C^{\one}_{i,j}$ and $C^{-\omega_{0}\one}_{i,j}$ 
for all $i,j$. Therefore by
Remark~\ref{rm:star-dual}, it follows that 
$t^{i-j}C^{\one}_{i,j} \in \mathcal{O}_{t}^{A_{0}}(SU(n+1))$
for all $1\leq i,j\leq n+1$. On the other hand, since
we also have $C^{\one}_{i,j} \in \mathcal{O}_{t}^{A_{0}}(SU(n+1))$,
it follows that
$t^{\min\{i-j,0\}}C^{\one}_{i,j}\in \mathcal{O}_{t}^{A_{0}}(SU(n+1))$
for all $1\leq i,j\leq n+1$.

Note that the matrix elements $C^{\one}_{i,j}$ satisfy the same
commutation relations as the elements $u_{i,j}$ in \cite{GirPal-2022tv}
(see equations~(2.6--2.9), \cite{GirPal-2022tv})  with
$q$ replaced by the variable $t$. Hence if one defines the
quantum determinant  (see Equation~(2.7), \cite{Koe-1991aa} or 
Chapter~4, \cite{ParWan-1991kr}) by
\begin{equation}\label{eq:comm-t5}
    D:=\sum_{\sigma\in\mathscr{S}_{n+1}}(-t)^{\ell(\sigma)}
    C^{\one}_{1,\sigma(1)}C^{\one}_{2,\sigma(2)}\dots C^{\one}_{n+1,\sigma(n+1)},
\end{equation}
with $\mathscr{S}_{n+1}$ being the permutation group on $n+1$ symbols and
$\ell(\sigma)$ being the length of the permutation $\sigma$, 
then one has $D=1$ and
\begin{align}\label{eq:c-onestar}
    (C^{\one}_{r,s})^{*} &= (-t)^{s-r}D^{r,s},
\end{align}
where $D^{s,r}$ denotes the $(s,r)$\raisebox{.2ex}{th} cofactor of the matrix
$(\!(C^{\one}_{i,j})\!)$.
Using the above equality, we will next prove that $(C^{\one}_{r,s})^{*}$ belongs
to the $A_{0}$-algebra generated by the elements
$t^{\min\{i-j,0\}}C^{\one}_{i,j}$, $1\leq i,j\leq n$. This will be done by
essentially reproducing the computations in the proof of Proposition~2.2.7
in \cite{GirPal-2022tv} in the present set up.

Define 
\[
i_{k}=\begin{cases}
          k & \text{if } 1\leq k < r,\\
          k+1 &\text{if } r\leq k \leq n,
      \end{cases} \qquad
j_{k}=\begin{cases}
          k & \text{if } 1\leq k < s,\\
          k+1 &\text{if } s\leq k \leq n,
      \end{cases}
\]
Let $\sigma\in\mathscr{S}_{n}$. Let
$\widetilde{\sigma}\in\mathscr{S}_{n+1}$ be the permutation that
takes $i_{k}\mapsto j_{\sigma(k)}$  for all $k\in\{1,2,\ldots,n\}$
and takes $r$ to $s$.
Then from (\ref{eq:c-onestar}), we have
\begin{align*}
  (C^{\one}_{r,s})^{*} &=
  (-t)^{s-r}\sum_{\sigma\in\mathscr{S}_{n}}(-t)^{\ell(\sigma)}
  C^{\one}_{i_{1},j_{\sigma(1)}}\ldots
  C^{\one}_{i_{n},j_{\sigma(n)}}\\
  &= (-t)^{s-r}\sum_{\sigma\in\mathscr{S}_{n}}
  (-t)^{\ell(\sigma)}C^{\one}_{1,\widetilde{\sigma}(1)}\ldots
  C^{\one}_{r-1,\widetilde{\sigma}(r-1)}
  C^{\one}_{r+1,\widetilde{\sigma}(r+1)}\ldots
  C^{\one}_{n+1,\widetilde{\sigma}(n+1)},
\end{align*}
For $r\leq s$, it follows immediately that $(C^{\one}_{r,s})^{*}$
belongs to the $A_{0}$-algebra generated by the elements
$t^{\min\{i-j,0\}}C^{\one}_{i,j}$, $1\leq i,j\leq n$.
Let us next assume that $r > s$. In this case, one has
{\allowdisplaybreaks
\begin{IEEEeqnarray*}{rCl}
    (C^{\one}_{r,s})^{*} &=& (-t)^{s-r}\sum_{\sigma\in\mathscr{S}_{n}}
  (-t)^{\ell(\sigma)}
  C^{\one}_{1,\widetilde{\sigma}(1)}\ldots 
  C^{\one}_{r-1,\widetilde{\sigma}(r-1)}
  C^{\one}_{r+1,\widetilde{\sigma}(r+1)}\ldots
   C^{\one}_{n+1,\widetilde{\sigma}(n+1)}\\
    &=& (-t)^{s-r}C^{\one}_{1,1}C^{\one}_{2,2}\ldots
  C^{\one}_{s-1,s-1}C^{\one}_{s,s+1}
  \ldots C^{\one}_{r-1,r}C^{\one}_{r+1,r+1}\ldots 
  C^{\one}_{n+1,n+1}\\
    & & {}+(-t)^{s-r}\sum_{\shortstack{$\scriptstyle \sigma\in\mathscr{S}_{n}$\\
   $\scriptstyle \sigma\neq id$}}
  (-t)^{\ell(\sigma)}C^{\one}_{1,\widetilde{\sigma}(1)}\ldots
  C^{\one}_{r-1,\widetilde{\sigma}(r-1)}
  C^{\one}_{r+1,\widetilde{\sigma}(r+1)}\ldots
   C^{\one}_{n+1,\widetilde{\sigma}(n+1)}\\
    &=& (-1)^{s-r}C^{\one}_{1,1}C^{\one}_{2,2}\ldots
  C^{\one}_{s-1,s-1}\left(t^{-1}C^{\one}_{s,s+1}\right)
  \ldots \left(t^{-1}C^{\one}_{r-1,r}\right)C^{\one}_{r+1,r+1}\ldots 
  C^{\one}_{n+1,n+1}\\
    & & {}+(-1)^{s-r}
t^{ s-r +{\sum_{k: j_{\sigma(k)}>i_{k}}}(j_{\sigma(k)}-i_{k})}\\
  && {}\hspace{2em}\times  
    \sum_{\shortstack{$\scriptstyle \sigma\in\mathscr{S}_{n}$\\
   $\scriptstyle \sigma\neq id$}}
  (-t)^{\ell(\sigma)}\left(t^{\min\{1-\widetilde{\sigma}(1),0\}}
        C^{\one}_{1,\widetilde{\sigma}(1)}\right)\ldots
  \left(t^{\min\{r-1-\widetilde{\sigma}(r-1),0\}}
       C^{\one}_{r-1,\widetilde{\sigma}(r-1)}\right) \\
 && {}\hspace{5em} \left(t^{\min\{r+1-\widetilde{\sigma}(r+1),0\}}
         C^{\one}_{r+1,\widetilde{\sigma}(r+1)}\right)\ldots
   \left(t^{\min\{n+1-\widetilde{\sigma}(n+1),0\}}
        C^{\one}_{n+1,\widetilde{\sigma}(n+1)}\right).
\end{IEEEeqnarray*}
}
By Lemma~2.2.6, \cite{GirPal-2022tv}, one has
$s-r +{\sum_{k: j_{\sigma(k)}>i_{k}}}(j_{\sigma(k)}-i_{k})\geq 0$
and therefore $(C^{\one}_{r,s})^{*}$ belongs to the 
$A_{0}$-algebra generated by the elements
$t^{\min\{i-j,0\}}C^{\one}_{i,j}$, $1\leq i,j\leq n$.
\eprf

We are now in a position to prove the main result of this section.
\bthm
\label{th:gpmy}
Let $C_{MY}$ denote the crystallization of the quantized function algebra
$C(SU_{q}(n+1))$ in the sense of Matassa \& Yuncken, and  $C_{GP}$
denote that in the sense of Giri \& Pal. Then one has $C_{MY}=C_{GP}$.
\ethm
\bprf
By Theorem~\ref{thm:MY-alternate}, $C_{MY}$ is the $C^{*}$-subalgebra of
$\mathcal{L}(\mathcal{H}_{Soi})$ generated by 
\[
\Bigl\{\lim_{q\to 0+}\psi^{(q)}_{Soi}\circ\vartheta_{q}^{}(a): a\in\otkazero\Bigr\}.
\]
By Theorem~\ref{th:otk}, it now follows that $C_{MY}$ is the $C^{*}$-subalgebra of
$\mathcal{L}(\mathcal{H}_{Soi})$ generated by
\[
\Bigl\{\lim_{q\to 0+}\psi^{(q)}_{Soi}\circ
  \vartheta_{q}^{}(t^{\min\{i-j,0\}}C^{\one}_{i,j}): 1\leq i,j\leq n\Bigr\}.
\]
From Theorem~4 and Remark~8.6 in \cite{GirPal-2024aa}, the above is precisely the
$C^{*}$-algebra $C_{GP}$.
\eprf


\noindent \textbf{Acknowledgement.} The second author would like to thank 
Prof.\ Sankaran Viswanath for making him aware of Lemma~\ref{lm:Connb}. 
He would also like to express his gratitude to Prof.\ Peter Littelmann 
for explaining his proof of the PRV conjecture in great detail.



\end{document}